\documentclass[11pt,twoside]{amsart}

\usepackage{amssymb}
\usepackage{amsmath}
\usepackage{epsfig}
\usepackage{bm}
\usepackage{mathrsfs}
\usepackage{bbm}
\usepackage{enumitem}

\usepackage[all]{xy}
\usepackage[dvipsnames]{xcolor}
\usepackage[colorlinks=true,citecolor=Plum,linkcolor=NavyBlue]{hyperref}
\usepackage{mathpazo}
\usepackage{empheq}

\newtheorem{theorem}{Theorem}[section] 

\newtheorem{proposition}[theorem]{Proposition}

\newtheorem{definition}[theorem]{Definition}

\def\build#1_#2^#3{\mathrel{\mathop{\kern 0pt#1}\limits_{#2}^{#3}}}

\newcommand*\mybox[1]{\colorbox{Gray!10}{\hspace{1em}#1\hspace{1em}}}

\def\A{{\mathscr A}}
\def\Ad{{\rm{Ad}}}
\def\C{{\mathbb C}}
\def\Conf{{\mathscr C}}
\def\d{{\; {\rm d}}}
\def\dw{{\rm d}}
\def\E{{\mathbb E}}
\def\F{{\mathbb F}}
\def\G{{\mathbb G}}
\def\g{{\mathfrak g}}
\def\geq{\geqslant}
\def\GL{{\mathrm{GL}}}
\def\hol{\text{hol}}
\def\id{{\rm id}}
\def\J{{\! \mathscr{J}}}

\def\leq{\leqslant}
\def\Loop{{\mathscr L}}
\def\Path{{\mathscr P}}
\def\R{{\mathbb R}}
\def\SO{{\mathrm{SO}}}
\def\Sp{{\mathrm{Sp}}}
\def\SU{{\mathrm{SU}}}
\def\Tr{{\rm Tr}}
\def\tr{{\rm tr}}
\def\U{{\mathrm U}}
\def\u{{\mathfrak u}}
\def\V{{\mathbb V}}
\def\YM{{\mathsf{YM}}}

\title[Two-dimensional quantum Yang--Mills theory and the Makeenko--Migdal equations]{Two-dimensional quantum Yang--Mills theory\\ and the Makeenko--Migdal equations}
\author{Thierry L\'evy}
\address{Laboratoire de Probabilit\'es, Statistique et Mod\'elisation (LPSM), Sorbonne Universit\'e}
\email{thierry.levy@sorbonne-universite.fr}
\date{\today}

\begin{document}

\maketitle

\setcounter{tocdepth}{2}
\tableofcontents

\section*{Introduction}

These notes, echoing a conference given at the Strasbourg-Zurich seminar in October 2017, are written to serve as an introduction to $2$-dimensional quantum Yang--Mills theory and to the results obtained in the last five to ten years about its so-called large $N$ limit. 

Quantum Yang--Mills theory, at least in the flavour that we will describe, combines differential geometric and probabilistic ideas. We would like to think, and hope to convince the reader, that this is less a complication than a source of beauty and enjoyment.

Some parts of our presentation will rely more distinctly on a probabilistic or a differential geometric background. We will however always try to keep technicalities aside and to favour explanation over demonstration. This is thus not, in the purest sense, a mathematical text: there will be essentially no proof. On the other hand, we will give fairly detailed examples of some computations that, we hope, are typical of the theory and illustrate it. 

Slightly different in aim and content, but also introductory, the notes \cite{LevySengupta} written with four hands with Ambar Sengupta can serve as counterpoint, or complement, to the present text. 

These notes are split in three parts. In the first, we explain the nature of the Yang--Mills holonomy process, which is the main object of interest of the theory. We do it from two perspectives, one differential geometric, and one probabilistic. This leads us to the definition of Wilson loop expectations, which are the most important numerical quantities of the theory. 

In the second part, we discuss several approaches to the computation of Wilson loop expectations, and illustrate them on several examples. The large $N$ limit of the theory makes a first appearance in this section, and we derive by hand some concrete instances of the Makeenko--Migdal equations which are the subject of the third part. We also included in the second part a discussion of the holonomy process on the sphere, and of the Douglas--Kazakov phase transition. 

In the third part, we describe the Makeenko--Migdal equations. In keeping with the style of these notes, we do not offer a proof of these equations, but we describe as carefully as we can Makeenko and Migdal's original derivation of them. Then, we discuss the amount of information carried by these equations and illustrate their power in the computation of the so-called master field, that is the large $N$ limit of Wilson loop functionals. \medskip

{\bf Acknowledgements.} I am grateful to Nalini Anantharaman and Ashkan Nikeghbali for organising the {\em 6th Strasbourg/Zurich - Meeting: Frontiers in Analysis and Probability}
and for their invitation to give the talk from which these notes are an expanded version. Part of the content of these notes was also covered in a series of three lectures that I gave in Lyon in June 2018 in a workshop on {\em Random matrices, maps and gauge theories} organised by Alice Guionnet, Adrien Kassel and Gr\'egory Miermont, whom I also want to thank. I am also  indebted to Adrien Kassel for his careful reading of a first version of this manuscript.

\section{Quantum Yang--Mills theory on compact surfaces}

\subsection{The holonomy process and the Yang--Mills action}\label{sec:HPYM}

The central object of study of quantum $2$-dimensional Yang--Mills theory is a collection of random unitary matrices indexed by the class $\Loop_{m}(M)$ of Lipschitz continuous loops based at some point $m$ on a compact surface $M$. This collection of random variables is called the {\em Yang--Mills holonomy process} and it is denoted by
\begin{equation}\label{eq:holproc}
(H_{\ell})_{\ell\in \Loop_{m}(M)}
\end{equation}
The idea of this collection of random variables arose, along a fairly convoluted path, from physical considerations relating to the description of certain kinds of fundamental interactions.\footnote{We will not describe this path, but indicate that it is marked by contributions of Chen Ning Yang and Robert Mills, the classical reference being \cite{YangMills}, of Alexander Migdal, who in \cite{Migdal} provided mathematicians with a usable description of a crucial part of Yang--Mills theory, of Leonard Gross who initiated a school of mathematical study of the $2$-dimensional Yang--Mills theory \cite{Gross1,Gross2,GKS}, of Bruce Driver and Ambar Sengupta, who finally gave in \cite{Driver,SenguptaAMS} the first mathematically rigorous definitions of the Yang--Mills holonomy process. This enumeration is of course much too short not to leave many important contributions aside: a more extensive bibliography can for instance be found in \cite{LevySengupta}.} It is, fortunately, not necessary to be familiar with the original motivation of Yang and Mills to understand what the Yang--Mills holonomy process is. 

In very broad terms, the basic data of the theory is a compact surface $M$ (for example a disk, a sphere, a cylinder, a torus) and a compact matrix group $G$ (for example $\U(1)$, $\SO(3)$, $\U(N)$). From this data, an infinite dimensional space of {\em connections} can be built\footnote{The exact nature of these connections can be ignored for the moment. If $G=\U(1)$, they can be pictured as magnetic potentials on $M$.}, on which an infinite dimensional symmetry group, the {\em gauge group} acts\footnote{In physical terms, two connections related by a gauge transformation represent two magnetic potentials corresponding to the same magnetic field.}, with infinite dimensional quotient, and one of the fundamental maps of the theory is the {\em holonomy map}
\[\xymatrix{
\{\text{connections}\}\big/\{\text{gauge group}\} \ar[rr]^{\text{\hspace{9mm}holonomy}} & &  {\rm Maps}(\Loop_{m}(M),G)\big/G }\]
On the right-hand side, the action of $G$ on the space of maps from $\Loop_{m}(M)$ to $G$ is by conjugation. Leaving this action aside, note that the distribution of the holonomy process \eqref{eq:holproc} is a probability measure on the space ${\rm Maps}(\Loop_{m}(M),G)$. We will call this space the space of {\em holonomies}.

One property that makes the holonomy map so important is that it is injective. It is thus legitimate to say that a connection is well described by its holonomy.

Another fundamental map of the theory is the {\em Yang--Mills action} $S_{\YM}$ which is a non-negative functional traditionally defined on the space of connections, but that can also be defined on the space of holonomies, so that the situation is
\begin{equation}\label{eq:diagS}
\xymatrix{
\{\text{connections}\}\big/\{\text{gauge group}\} \ar[rr]^{\text{\hspace{9mm}holonomy}} \ar[dr]_{S_{\YM}} & &  {\rm Maps}(\Loop_{m}(M),G)\big/G \ar[dl]^{S_{\YM}}\\
& [0,\infty] &
}
\end{equation}
The {\em Yang--Mills measure} is heuristically described as the Boltzmann probability measure, on the space of connections or on the space of holonomies, associated with the Yang--Mills action. The typical formula that one finds in the literature is 
\begin{equation}\label{eq:YM}
\d\mu_{\YM}(\omega)=\frac{1}{Z} \, e^{-\frac{1}{2T}S_{\YM}(\omega)}\, \d\omega
\end{equation}
where $T$ is a positive real parameter called the coupling constant. Here, $\omega$ is meant to stand for a connection or for a holonomy, depending on one's preferred point of view. This expression is however plagued with difficulties:  on the infinite dimensional spaces where the Yang--Mills measure is supposed to live, there is no Lebesgue-like reference measure that could reasonably play the role of $\d\omega$, and even if there were, one would not expect the Yang--Mills measure to be absolutely continuous with respect to it; moreover, because of the action of the gauge group, the most sensible value for the normalisation constant would be $Z=+\infty$; and one does finally not expect a typical $\omega$ in the sense of the Yang--Mills measure to be regular enough to have a finite Yang--Mills action. 

One of the goals of the $2$-dimensional quantum Yang--Mills theory is to find a way of sorting out these difficulties and to construct rigorously a probability measure that can honestly be called the Yang--Mills measure. The situation may look rather desperate, but it is uplifting to realise that after replacing the space of connections, or holonomies, by a space of real-valued functions on $[0,1]$ and the Yang--Mills action by the square of the Sobolev $H^{1}$ norm, the analogous problem is almost just as ill-posed but has a very well-known solution, namely the Wiener measure. The main difference between the Wiener and the Yang--Mills cases is the presence in the latter of the gauge symmetry. Symmetry can however be a nuisance or a guide, and it turns out to be possible, in Yang--Mills theory, to make gauge symmetry an ally rather than a foe. 

We will now describe more precisely the three maps appearing in the diagram \eqref{eq:diagS}. The holonomy map and the Yang--Mills action on the space of connections are differential geometric in nature. We start by describing them, and then turn to the Yang--Mills action on the space of holonomies. It would be unfair to say that the content of Section \ref{sec:YMconn} can safely be completely ignored: we will refer to it later, in particular in Section \ref{sec:MMproof}. However, it is certainly possible to skip it at first reading and to jump to Section \ref{sec:YMholo}. 

\subsection{The Yang--Mills action: connections}\label{sec:YMconn}

In this section, we assume from the reader some familiarity with the differential geometry of principal bundles. We give brief reminders of the main definitions, but this is of course not the place for a complete exposition. For details, and although some might find it too Bourbakist in style, we recommend the second chapter of the first volume of the classical opus by Kobayashi and Nomizu \cite{KobayashiNomizu}.

\subsubsection{The Yang--Mills action}\label{sec:YMA}

Although we are concerned in this text with compact surfaces, we will describe the Yang--Mills action in the more general context of compact Riemannian manifolds of arbitrary dimension --- this is not more difficult.

Let $M$ be a compact connected Riemannian manifold. Let $G$ be a compact Lie group with Lie algebra~$\g$. Assume that $\g$ is endowed with a scalar product $\langle \cdot,\cdot \rangle$ that is invariant under the adjoint representation $\Ad:G\to \GL(\g)$.\footnote{\label{fn:GU}The typical example that we have in mind is $G=\U(N)$ and, for all $X,Y\in \u(N)$ skew-Hermitian $N\times N$ matrices, $\langle X,Y\rangle = N \Tr(X^{*}Y)$.}
Let $\pi:P\to M$ be a principal $G$-bundle over $M$.\footnote{The manifold $P$ is thus acted on, on the right, by $G$. For small open subsets $U$ of $M$, the part $\pi^{-1}(U)$ of the manifold $P$ that sits above $U$ is equivariantly diffeomorphic to $U\times G$, with $\pi$ being the first coordinate map and $G$ acting by translations on the right on the second coordinate. A principal bundle is {\em trivial} if it is globally isomorphic to $M\times G$.} Let $\A$ denote the space of connections on $P$. It is an affine subspace of the space of $\g$-valued differential $1$-forms on $P$. For every connection $\omega\in \A$, the curvature of $\omega$ is the form $\Omega=d\omega+\frac{1}{2}[\omega\wedge \omega]$.\footnote{This definition of the curvature is made slightly ambiguous by the coexistence, in the literature, of two different conventions regarding the definition of the exterior product and the exterior differential of differential forms. Since it took me some time to clarify this elementary point, I want to record it here, to the price of a  rather long footnote.

The two conventions could be called `simplicial' and `cubical' according to their respective definitions of the exterior product of $1$-forms:
\[(\alpha_{1}\wedge \ldots \wedge \alpha_{k})(X_{1},\ldots,X_{k})=\left\{\begin{array}{rl} \frac{1}{k!}\det\big[(\alpha_{i}(X_{j}))_{i,j=1\ldots k}\big] & \text{(simplicial)}\\
\det\big[(\alpha_{i}(X_{j}))_{i,j=1\ldots k}\big] & \text{(cubical)}
\end{array}\right.\]
Each convention is supported by illustrious authors, including, for the simplicial one, Kobayashi and Nomizu \cite[p. 35]{KobayashiNomizu} and Morita \cite[Eq. (2.14) p. 70]{Morita}, and for the cubical one, Spivak \cite[p. 203]{Spivak}. Since everyone agrees on the formula $d(\alpha\wedge \beta)=d\alpha \wedge \beta+(-1)^{{\rm deg}(\alpha){\rm deg}(\beta)} \alpha \wedge d\beta$, there must also be two competing definitions of the exterior differential. Specifically, the two definitions are related by the formula $d^{\text{simplicial}}\alpha=\frac{1}{{\rm deg}(\alpha)+1}d^{\text{cubical}}\alpha$ (compare, for instance, \cite[p. 36]{KobayashiNomizu} or \cite[Thm. 2.9 p. 71]{Morita} and \cite[Thm 13 p. 213]{Spivak}). The formula $d\alpha(X,Y)=X\alpha(Y)-Y\alpha(X)-\alpha([X,Y])$, for instance, belongs to the cubical school.

Returning to the definition of the curvature, it has a different meaning with each convention, but fortunately, the simple relation $\Omega^{\text{simplicial}}=\frac{1}{2}\Omega^{\text{cubical}}$ holds. Let us be more explicit about this definition: the expression $\omega\wedge \omega$ is to be understood as a $\g\otimes \g$-valued $2$-form, which is then composed by the Lie bracket to yield a $\g$-valued $2$-form. Explicitly, if $X$ and $Y$ are two vector fields defined on an open subset of $P$, then the curvature of $\omega$ is defined on this open set by
\[\Omega^{\text{cubical}}(X,Y)=2\Omega^{\text{simplicial}}(X,Y)=X\omega(Y)-Y\omega(X)-\omega([X,Y])+[\omega(X),\omega(Y)]\]
Note that there is universal agreement on what it means for the curvature to vanish.

Finally, since everyone also agrees that Stokes' formula is free of any coefficient, each convention on the definition of the exterior differential entails its own definition of the integral. This is slightly hidden by the fact that everyone agrees on the formula $\int_{[0,1]^{n}} dx_{1}\wedge \ldots \wedge dx_{n}=1$ (see \cite[Sec. 3.2 (a), p. 104]{Morita} and \cite[Prop. 1 p. 247]{Spivak}), but it must be realised that the differential form that is denoted by $dx_{1}\wedge \ldots \wedge dx_{n}$ is not the same for everyone. Specifically, the relation is $\int^{\text{simplicial}} \alpha = {\rm deg}(\alpha)! \int^{\text{cubical}}\alpha$.

Finally, there is agreement on the meaning of the curvature as a linear map from the space of smooth $2$-chains in~$P$ to $\g$. 
}
 This~$\g$-valued $2$-form on $P$ vanishes on vertical vectors and is $G$-equivariant. It can thus be seen as a~$2$-form on $M$ with values in the adjoint bundle $\Ad (P)$. Using the Hodge operator of the Riemannian structure of $M$, one can form the $(\Ad (P)\otimes \Ad (P))$-valued form of top degree~$\Omega \wedge {\star}\Omega$ on $M$. Contracting this form with the Euclidean structure of $\Ad(P)$ induced by the invariant scalar product on $\g$ yields the real-valued differential form of top degree $\langle \Omega \wedge {\star}\Omega\rangle$. This form can be integrated\footnote{The definition of the Yang--Mills action seems to require an orientation of $M$. In fact, this orientation is used twice, once to define the Hodge dual $\star \Omega$ of $\Omega$ and once to integrate $\langle \Omega,\star \Omega\rangle$ over $M$. Reversing the orientation changes the Hodge dual and the integral by a sign,
so that if $M$ is orientable, the definition of $S_{\YM}$ is independent of the choice of orientation of $M$. Moreover, if $M$ is not orientable, $S_{\YM}$ can still be defined using a partition of unity.} to yield the Yang--Mills action of $\omega$:
\begin{equation}\label{eq:defSint}
S_{\YM}(\omega)=\frac{1}{2}\int_{M} \langle \Omega \wedge {\star}\Omega\rangle
\end{equation}
In words, the Yang--Mills action of a connection is nothing more than one half of the squared~$L^{2}$ norm of its curvature.\footnote{Considering that the curvature is a kind of derivative of the connection, the Yang--Mills action stands thus in close analogy with the squared $H^{1}$ norm of a real-valued function on $[0,1]$.}

Let us describe locally, in coordinates, the scalar function that is integrated over $M$ to compute $S_{\YM}(\omega)$. For this, let us consider an open subset $U$ of $M$ on which there exist local coordinates $x_{1},\ldots,x_{n}$ on $M$ and over which $P$ is trivial. Let us choose a section\footnote{To say that $\sigma$ is a section of $P$ over $U$ means  that $\pi\circ \sigma=\id_{U}$. The existence of such a section is equivalent to the triviality of the restriction of $P$ over $U$. In particular, the existence of a global section $\sigma:M\to P$ is equivalent to the triviality of the bundle $\pi:P\to M$. The reader who is more familiar with vector bundles than principal bundles might at first be surprised by this statement, since a vector bundle can admit a global section, even a non-vanishing one, without being trivial. However, the existence of a section for a principal bundle corresponds, for a vector bundle, to the existence of a basis of sections.
} $\sigma:U\to P$ of~$P$ over $U$. Let us define $A=\sigma^{*}\omega$. Then in the local coordinates on $U$, the $1$-form $A$ writes $A_{1}\d x_{1}+\ldots+A_{n}\d x_{n}$, where $A_{1},\ldots,A_{n}$ are maps from $U$ to $\g$. Then $F=\sigma^{*}\Omega$ writes
\[F=\sum_{1\leq i<j\leq n} \big(\partial_{i}A_{j}-\partial_{j}A_{i}+[A_{i},A_{j}]\big)\d x_{i}\wedge \dw x_{j}\]
and the contribution of $U$ to the Yang--Mills action of $\omega$ is
\[\frac{1}{2}\int_{U}\langle \Omega \wedge {\star}\Omega\rangle=\frac{1}{2}\sum_{1\leq i<j\leq n} \int_{U} \big\| \partial_{i}A_{j}-\partial_{j}A_{i}+[A_{i},A_{j}]\big\|^{2} \d\text{vol}(x)\]
where $\d\text{vol}(x)$ is the Riemannian volume measure on $M$, and $\|\cdot\|$ is the Euclidean norm on~$\g$ associated with the invariant scalar product $\langle \cdot,\cdot \rangle$. The analogy with the squared Sobolev $H^{1}$ norm should be even more obvious on this expression.

\subsubsection{Gauge transformations}\label{sec:gaugetransfo}

The gauge group, that we denote by $\J$, is the group of $G$-equivariant diffeomorphisms of $P$ over the identity of $M$.\footnote{An element $j$ of the gauge group is a diffeomorphism $j:P\to P$ that leaves each fibre of $P$ globally stable, and acts on it in a way that commutes with the action of $G$ on the right on $P$. For the bundle $P=M\times G\to M$, the gauge group can be identified with $\J=C^{\infty}(M,G)$ acting pointwise on $P$ by multiplication on the left on the second coordinate.} It acts by pull-back on $\A$ and a routine verification shows that it leaves $S_{\YM}$ invariant. Thus, the Yang--Mills action descends to a function
\[S_{\YM}:\A\big/\J \to [0,\infty)\]
the study of which is the subject of classical Yang--Mills theory. 

Let us display the formulas which give, through a local section of $P$, the action of the gauge group on a connection and its curvature. These formulas are indeed useful, and ubiquitous in the literature. Let $j:P\to P$ be a gauge transformation. Let $\sigma:U\to P$ be a local section of~$P$ over an open subset $U$ of $M$. Then there exists a unique function $g:U\to G$ such that for every $x\in U$, one has $j(\sigma(x))=\sigma(x)g(x)$. Then, letting $j$ act on a connection $\omega$ yields the new connection $j\cdot \omega=j^{*}\omega$ and transforms on one hand $A$ into
\[A^{g}=\sigma^{*}(j\cdot\omega)=g^{-1}Ag+g^{-1}\d g\]
and on the other hand $F$ into
\[F^{g}=g^{-1}Fg\]
This formula explains the invariance of the Yang--Mills action: without trying to be perfectly precise, one can say that the action of a gauge transformation conjugates the curvature at each point of $M$ by some element of $G$, and thus leaves its Euclidean norm unchanged.

\subsubsection{Some questions of classical Yang--Mills theory}

Let us mention, without giving any details, a few examples of the questions that arise in the study of the Yang--Mills action.

\begin{itemize}[itemsep=0pt]
\item The set $S_{\YM}^{-1}(0)$ is the {\em moduli space of flat connections}, that is, the quotient of the set of flat connections by the action of the gauge group. It is a finite-dimensional orbifold with a rich geometric structure, the study of which is both an old and an active subject of investigation \cite{Goldman1,Goldman2,Witten1,Witten2,KingSengupta,Liu1,Liu2}. 
\item The Yang--Mills action can be understood as arising, through appropriate reformulation and generalisation, from a Lagrangian formulation of Maxwell's equations of the electromagnetic field. The critical points of the Yang--Mills action are thus of special interest: they are, in a sense, the classical physical fields of Yang--Mills theory. They are called {\em Yang--Mills connections} and a milestone in their study in the $2$-dimensional case is \cite{AtiyahBott}.
\item When $M$ is $4$-dimensional, the Yang--Mills action is conformally invariant, in the sense that it depends on the Riemannian metric on $M$ only through its conformal class. There is an extensive literature devoted to Yang--Mills connections on $4$-dimensional manifolds \cite{HitchinBourbaki}. Looking for self-dual Yang--Mills connections on $\R^{4}$ that are invariant by translation in two directions, for example, leads to the study of Hitchin equations and Higgs bundles \cite{Hitchin}. 
\item From a physical point of view, the Yang--Mills action of a connection is an appropriate measure of its non-triviality. From an analytical point of view, however, it turns out that a natural way of measuring a connection is its Sobolev $H^{1}$ norm.\footnote{Here, we are talking about connections as elements of $\A$, not of the quotient $\A/\J$.} The Yang--Mills action is controlled by the $H^{1}$ norm, but not conversely. A flat connection, that is, a connection with Yang--Mills action $0$, can be given an arbitrarily large $H^{1}$ norm by an appropriate gauge transformation. A beautiful theorem of Karen Uhlenbeck states that level sets of the Yang--Mills action, that is, the sets of the form $\{S_{\YM}\leq c\}$, $c\in \R_{+}$, are sequentially weakly compact in $H^{1}$ up to gauge transformation: from any sequence of connections with bounded Yang--Mills action, one can extract a subsequence which, after suitable gauge transformation of each term, converges weakly in $H^{1}$ \cite{Uhlenbeck}.
\item The Yang--Mills action gives rise to a gradient flow, which formally is the solution of the differential equation $\partial_{t}\omega_{t}=-\nabla_{\omega_{t}}S_{\YM}$. This is the {\em Yang--Mills flow} \cite{Rade}. There is currently an active investigation of stochastic perturbations of this flow in cases where~$M$ is $2$- or $3$-dimensional \cite{Shen, Chevyrev}. 
\end{itemize}

\subsubsection{The holonomy map}\label{subsec:holomap}

A fundamental construction associated with a connection is that of the {\em holonomy}, or {\em parallel transport}, that it induces. For every continuous and piecewise smooth curve $c:[0,1]\to M$, the parallel transport along $c$ determined by the connection $\omega$ is the $G$-equivariant mapping $\hol(\omega,c):P_{c_{0}}\to P_{c_{1}}$ which to every point $p$ of $P_{c_{0}}$ associates the endpoint of the unique continuous curve $\tilde c:[0,1]\to P$ such that $\tilde c_{0}=p$, $\pi\circ \tilde c=c$ and for all $t\in [0,1]$ at which $c$ is differentiable, $\omega(\dot{\tilde{c}}_{t})=0$.

This parallel transport enjoys the following properties, which are of fundamental importance.
\begin{itemize}[itemsep=0pt]
\item It is unaffected by a change of parametrisation of the curve.
\item If $c:[0,1]\to M$ is a curve and $c^{-1}$ denotes the same curve traced backwards, that is, $c^{-1}_{t}=c_{1-t}$, then $\hol(\omega,c^{-1})=\hol(\omega,c)^{-1}$.
\item If $c$ and $c'$ are two curves such that $c_{1}=c'_{0}$, so that the concatenation $cc'$ is well defined, then $\hol(\omega,cc')=\hol(\omega,c')\circ \hol(\omega,c)$.
\end{itemize}

It will be useful to understand a bit more concretely how this parallel transport can be computed, and how it gives rise to a holonomy in the sense that we gave to this word in Section~\ref{sec:HPYM}. 

\begin{figure}[h!]
\begin{center}
\includegraphics{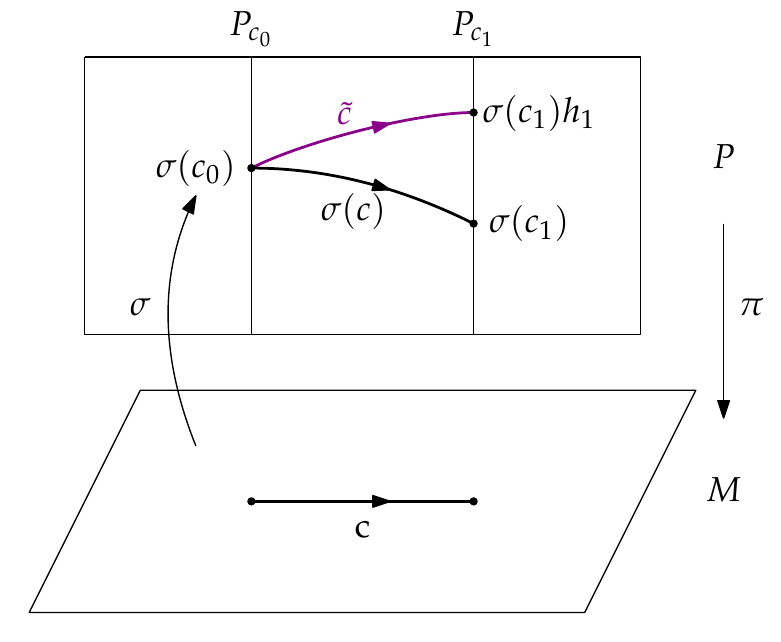}
\caption{\small \label{fig:holonomie} The difference between the horizontal lift of $c$ starting at $\sigma(c_{0})$, denoted in this picture by $\tilde c$, and $\sigma(c)$, the image of $c$ by the local section $\sigma$, is measured by the function $h$ which solves \eqref{eq:defh}.}
\end{center}
\end{figure}

Assume that the range of the curve $c$ lies in an open subset $U$ of $M$ over which the fibre bundle $P$ is trivial.\footnote{If $c$ does not lie in such an open subset, it can be split into finitely many pieces which do and the holonomy along $c$ is simply the product of the holonomies along these shorter pieces.} Let $\sigma:U\to P$ be a section of $P$ over $U$. Set $A=\sigma^{*}\omega$. It is a $1$-form on $U$ with values in $\g$. The solution of the differential equation
\begin{equation}\label{eq:defh}
\dot h_{t}=-A(\dot c_{t})h_{t}, \ h_{0}=1_{G}
\end{equation}
is a curve $h:[0,1]\to G$ which starts from the unit element $1_{G}$. The endpoint of this curve computes the parallel transport along $c$ determined by $\omega$ in the sense that
\[\hol(\omega,c)(\sigma(c_{0}))=\sigma(c_{1})h_{1}\]
This relation is illustrated in Fig. \ref{fig:holonomie}.

Let us introduce the notation
\[\hol_{\sigma}(\omega,c)=h_{1}\] 
the holonomy of $\omega$ along $c$ read in the section $\sigma$. This object has the drawback of depending on the choice of a local section of the bundle, but the great advantage of being fairly concrete, namely an element of $G$, that is, in many situations, a matrix.

If $j\in \J$ is a gauge transformation of $P$, recall from Section \ref{sec:gaugetransfo} that $j\cdot \omega=j^{*}\omega$ is the pull-back of $\omega$ by the diffeomorphism $j$ of $P$. Then the holonomy of $j\cdot\omega$ along $c$ is related to that of $\omega$ by the relation 
\[\hol(j\cdot \omega,c)=j_{|P_{c_{1}}}^{-1}\circ \hol(\omega,c)\circ j_{|P_{c_{0}}}\]
Through the local section $\sigma:U\to M$, and letting $g:U\to G$ be the function such that $j(\sigma(x))=\sigma(x)g(x)$ for every $x\in U$, this relation takes the more explicit form
\begin{equation}\label{eq:gaugehol}
\hol_{\sigma}(j\cdot \omega,c)=g_{c_{1}}^{-1}\hol_{\sigma}(\omega,c)g_{c_{0}}
\end{equation}

It follows from \eqref{eq:gaugehol} that for all loop $\ell$ on $M$, that is, all curve which ends at its starting point, the conjugacy class of $\hol_{\sigma}(\omega,\ell)$ is not affected\footnote{Incidentally, this class does not depend on the local section $\sigma$ either. 
} by a gauge transformation of $\omega$. 

More generally, given a base point $m$ on $M$, and denoting by $\Loop_{m}^{\infty}(M)$ the class of piecewise smooth loops on $M$ based at $m$, the orbit of 
\[(\hol_{\sigma}(\omega,\ell) : \ell\in \Loop_{m}^{\infty}(M))\in {\rm Maps}(\Loop_{m}^{\infty}(M),G)\]
under the action of $G$ by simultaneous conjugation is not affected by a gauge transformation of~$\omega$. This explains how a connection modulo gauge transformation defines a holonomy modulo conjugation. 

The following result makes precise the statement that the horizontal arrow of \eqref{eq:diagS} is injective.

\begin{theorem}\label{thm:carconnJ} Let $m$ be a point of $M$. Let $\sigma$ be a section of $P$ in a neighbourhood of $m$. For any two connections $\omega$ and $\omega'$ on $P$, the following assertions are equivalent.\\
\indent 1. There exists a gauge transformation $j\in \J$ such that $j\cdot \omega=\omega'$.\\
\indent 2. There exists $g\in G$ such that for all loop $\ell\in \Loop_{m}^{\infty}(M)$, the equality $\hol_{\sigma}(\omega',l)=g^{-1}\hol_{\sigma}(\omega,l)g$ holds.
\end{theorem}

\subsection{The Yang--Mills action: holonomies} \label{sec:YMholo}

We will now give an alternative of the Yang--Mills action that is less classical and, most importantly, specific to the $2$-dimensional case. To give an idea of the nature of this second description, let us pursue the analogy with the Wiener measure and the Sobolev $H^{1}$ norm. Consider a smooth function $b:[0,1]\to \R$ with $b(0)=0$. The squared $H^{1}$ norm of $b$ can be expressed at least in the following two ways:
\begin{equation}\label{eq:H1}
\|b\|_{H^{1}}^{2} =\int_{0}^{1} |\dot b(t)|^{2} \d t=\sup_{0\leq t_{0}< t_{1}<\ldots<t_{n}\leq 1} \sum_{k=1}^{n} \frac{|b(t_{k})-b(t_{k-1})|^{2}}{t_{k}-t_{k-1}}
\end{equation}
The integral expression corresponds to the description of the Yang--Mills action that we gave in the last section and is very similar to \eqref{eq:defSint}. We will now give another description, similar to the second, more combinatorial one.
\subsubsection{Holonomies}

The main algebraic property of the holonomy of a connection, already mentioned in Section \ref{subsec:holomap}, is that it is a {\em multiplicative} map from $\Loop_{m}^{\infty}(M)$ to $G$. Let us formulate this in a slightly different way. 

Recall that $M$ is a compact Riemannian manifold and $G$ a compact Lie group. Let $\Path(M)$ denote the set of all Lipschitz continuous\footnote{In this text, we consider alternatively paths that are piecewise smooth and paths that are Lipschitz continuous. We do so for reasons of technical convenience, and the reader should not be overly worried by what can safely be regarded as a secondary issue.} paths on $M$, two paths being identified if they differ only by an increasing change of parametrisation. Let us call a function $h:\Path(M)\to G$ {\em multiplicative} if it satisfies the following two properties.
\begin{itemize}[itemsep=0pt]
\item For all path $c$, letting $c^{-1}$ denote the same path traced backwards, one has $h(c^{-1})=h(c)^{-1}$.
\item For all paths $c$ and $c'$ such that $c$ finishes where $c'$ starts, so that the concatenated path $cc'$ is defined, one has $h(cc')=h(c')h(c)$.
\end{itemize}
More generally, given a subset $P$ of $\Path(M)$, we say that a function $h:P\to G$ is multiplicative if it satisfies the above conditions whenever all the paths involved belong to the subset $P$.

Let us denote by ${\rm Mult}(\Path(M),G)$ (resp. by ${\rm Mult}(P,G)$) the subset of ${\rm Maps}(\Path(M),G)$ (resp. of ${\rm Maps}(P,G)$) formed by all multiplicative maps. 

There is an action of the {\em gauge group} ${\rm Maps}(M,G)$ on ${\rm Mult(\Path(M),G)}$ defined as follows. Consider $g:M\to G$ and a multiplicative map $h:\Path(M)\to G$. For all path $c$ starting at $c_{0}$ and finishing at $c_{1}$, define
\begin{equation}\label{eq:gaugelattice}
(g\cdot h)(c)=g_{c_{1}}^{-1}h(c)g_{c_{0}}
\end{equation}
an equation that should be compared with \eqref{eq:gaugehol}. It is not difficult to check that the map $g\cdot h$ is still multiplicative.

Let $m$ be a point of $M$. A multiplicative function can be restricted to $\Loop_{m}(M)$ and the action of ${\rm Maps}(M,G)$ on this restricted map reduces to the action of $G$ by conjugation. The following fact may seem surprising at first sight, but it is not difficult to prove.

\begin{proposition}\label{prop:multmult} For all $m\in M$, the restriction map 
\[{\rm Mult}(\Path(M),G)\big/ {\rm Maps}(M,G) \longrightarrow {\rm Mult}(\Loop_{m}(M),G)\big/ G\]
is a bijection.
\end{proposition}

We call either side of this bijection the space of {\em holonomies}. Thanks to the multiplicativity and the gauge symmetry, a holonomy can either be seen as a group-valued function on the set of all paths, or on the set of all loops based at some reference point $m$ on $M$.

\subsubsection{Graphs on surfaces}\label{sec:graphs}

We will now assume that $M$ is a $2$-dimensional manifold: it is thus a compact surface. We announced an expression of the Yang--Mills action similar to the rightmost term of \eqref{eq:H1}: the role of subdivisions of the interval $[0,1]$ will be played by {\em graphs} on $M$. This will be the occasion of a first encounter with this notion that is central to the construction of the~$2$-dimensional Yang--Mills measure.

Let us call {\em edge} an element of $\Path(M)$ that is injective --- note that this does not depend on the way in which the path is parametrised. A {\em graph} is a finite set of edges, stable by the reversal map $e\mapsto e^{-1}$, and in which any two edges either form a pair $\{e,e^{-1}\}$, or meet, if at all, at some of their endpoints. 

The {\em vertices} of a graph are the endpoints of its edges. The {\em faces} of a graph are the connected components of the complement in $M$ of the union of its edges. A graph is conveniently described as a triple $\G=(\V,\E,\F)$ consisting of a set of vertices, a set of edges and a set of faces, but it is in fact entirely determined by the set $\E$ of its edges. 

A crucial additional assumption is that every face of a graph must be homeomorphic to a disk. This guarantees that the $1$-skeleton of the graph correctly represents the topology of the surface, to the extent that a $1$-dimensional object can represent a $2$-dimensional one.

\subsubsection{The Yang--Mills action}\label{sec:YMac}

Let $\G$ be a graph on our compact surface $M$. We will denote by$\Path(\G)$ the set of paths that can be constructed as concatenations of edges of $\G$. To each face $F$ of $\G$, we can associate in an almost unequivocal way a loop $\partial F$ that winds exactly once around$F$. To give a perfectly rigourous definition of this loop is less simple than one might expect, but there is nothing counterintuitive in it. It is only {\em almost} well defined because there is no preferred starting point for this loop. However, if $f:\Path(\G)\to G$ is a multiplicative function, then the conjugacy class of the element $h(\partial F)$ of $G$ is well defined. In particular, the Riemannian distance, in $G$, between the element $h(\partial F)$ and the unit element $1_{G}$, is well defined.\footnote{This distance is defined by the bi-invariant Riemannian metric on $G$ associated with the invariant scalar product chosen on its Lie algebra, see the first lines of Section \ref{sec:YMA}.} This distance is, moreover, not affected  by the action of an element of the gauge group ${\rm Maps}(M,G)$ on $h$.

We can now define the Yang--Mills action on the space of holonomies by setting, for all $h\in {\rm Mult}(\Path(M),G)$,
\begin{equation}\label{eq:YMH}
S_{\YM}(h)=\sup\bigg\{\sum_{F\in \F} \frac{d_{G}(1_{G},h(\partial F))^{2}}{\text{area}(F)}\  :\  \G \text{ graph on } M\bigg\} 
\end{equation}
where the area of a face $F$ is computed using the Riemannian structure on $M$.

It is manifest on this expression that, in the case where $M$ is a surface, the only part of the Riemannian structure on $M$ that is used in the definition of the Yang--Mills action is the Riemannian volume, in this case the Riemannian area. This is of course also true, be it in a slightly less apparent way, of the definition \eqref{eq:defSint}. 

\begin{proposition} Assume that $M$ is $2$-dimensional. Then the definitions \eqref{eq:defSint} and \eqref{eq:YMH} of the Yang--Mills action agree. More precisely, for every connection $\omega$ inducing a holonomy $h$, the equality $S_{\YM}(\omega)=S_{\YM}(h)$ holds.
\end{proposition}

\subsection{The Yang--Mills holonomy process}

We will now explain how to construct the Yang--Mills holonomy process. Although the definition of this process is derived, at a heuristic level, from the Yang--Mills action, the process and the action are logically unrelated. We can thus start afresh, from a compact surface $M$ on which we have a Riemannian structure, or at least a measure of area, and a compact Lie group $G$, on the Lie algebra of which we have an invariant scalar product. 

\subsubsection{The configuration space of lattice Yang--Mills theory}

One piece of information that we need to retain from the previous sections is the notion of graph on our surface $M$ (see Section \ref{sec:graphs}). Let us choose a graph $\G=(\V,\E,\F)$ on $M$.
The {\em configuration space} associated to a graph $\G$ on our surface $M$ is the manifold
\[\Conf_{\G}=\{g=(g_{e})_{e\in \E} \in G^{\E} : \forall g\in G, g_{e^{-1}}=g_{e}^{-1}\}={\rm Mult}(\E,G)\]
of all ways of assigning an element of $G$ to each oriented edge, in a way that is consistent with the orientation reversal.

Recall that we denote by $\Path(\G)$ the set of paths that can be constructed as concatenations of edges of $\G$. The configuration space $\Conf_{\G}$ is naturally in one-to-one correspondence with the set ${\rm Mult}(\Path(\G),G)$ of all multiplicative maps from $\Path(\G)$ to $G$.

Choosing an {\em orientation} of $\G$, that is, a subset $\E^{+}\subset \E$ containing exactly one element in each pair $\{e,e^{-1}\}$ allows one to realise the configuration space in the slightly less canonical, but easier to handle, way
\[\Conf_{\G}=G^{\E^{+}}\]
This makes it easy, for instance, to endow $\Conf_{\G}$ with a probability measure, namely the Haar measure on $G^{\E^{+}}$. The invariance of the Haar measure on the compact group $G$ under the inverse map $x\mapsto x^{-1}$ implies that this measure on $\Conf_{\G}$ does not depend on the choice of orientation. We denote it by $\dw g$.

Every path $c\in \Path(\G)$ can be uniquely written as a concatenation of edges $c=e_{1}^{\epsilon_{1}}\ldots e_{n}^{\epsilon_{n}}$ with $e_{1},\ldots,e_{n}\in \E^{+}$ and $\epsilon_{1},\ldots,\epsilon_{n}\in \{-1,1\}$. To such a path $c=e_{1}^{\epsilon_{1}}\ldots e_{n}^{\epsilon_{n}}$ we associate a holonomy map
\begin{align}\label{eq:defholo}
h_{c}: \Conf_{\G} & \longrightarrow G \\ \nonumber
g & \longmapsto g_{e_{n}}^{\epsilon_{n}}\ldots g_{e_{1}}^{\epsilon_{1}}
\end{align}
Our goal is to endow the configuration space $\Conf_{\G}$ with an interesting probability measure, so as to make the collection of maps $(h_{c})_{c\in \Path(\G)}$ into a collection of $G$-valued random variables.

\subsubsection{The Driver--Sengupta formula} \label{sec:DS}

In order to define this probability measure, we need to introduce the heat kernel on $G$, or more accurately the fundamental solution of the heat equation. The invariant scalar product on the Lie algebra $\g$ determines a bi-invariant Riemannian structure on $G$, and a Laplace-Beltrami operator $\Delta$. We consider the function $p:\R^{*}_{+}\times G \to \R^{*}_{+}$ that is the unique positive solution of the heat equation $(\partial_{t}-\frac{1}{2} \Delta)p=0$ with initial condition $p(t,x)\d x\Rightarrow \delta_{1_{G}}$ as $t\to 0$. We use the notation $p_{t}(x)=p(t,x)$. A crucial property of this function is that, for all $t>0$ and all $x,y\in G$, we have $p_{t}(yxy^{-1})=p_{t}(x)$. We refer to this property as the {\em invariance under conjugation} of the heat kernel.

We mentioned at the end of Section \ref{sec:YMac} that, in the $2$-dimensional setting, the Yang--Mills action depends on a Riemannian structure of the surface $M$ only through the Riemannian area that it induces. We will denote by $|F|$ the area of a Borel subset $F$ of $M$. 

Given a face $F$ of our graph, recall that we denote by $\partial F$ a path that goes once around this face in the positive direction. Recall also that this path is ill-defined because there is no preferred vertex on the boundary of $F$ from which to start it. However, this indeterminacy only results in an indeterminacy {\em up to conjugation} for the holonomy map $h_{\partial F}$. Thanks to the invariance under conjugation of the heat kernel, the function $g\mapsto p_{t}(h_{\partial F}(g))$ is still well defined on $\Conf_{\G}$ for every~$t>0$.

We can now write the formula which is the basis of the definition of the $2$-dimensional Yang--Mills measure. It is due to Bruce Driver in the case where $M$ is the plane, or a disk, and to Ambar Sengupta when $M$ is an arbitrary compact surface. Recall that $T$ is a positive real parameter of the measure. We define, on $\Conf_{\G}$, the probability measure
\begin{empheq}[box=\mybox]{equation}
\label{eq:discreteYM}\tag{DS}
\dw\mu^{\G,T}_{\YM}(g)=\frac{1}{Z(\G,T)} \prod_{F\in \F} p_{T|F|}(h_{\partial F}(g))\d g
\end{empheq}
Here, $Z(\G,T)$ is the normalisation constant that makes $\mu^{\G,T}_{\YM}$ a probability measure on $\Conf_{\G}$.

The gauge group ${\rm Maps}(\V,G)$ acts on the configuration space $\Conf_{\G}$ by a formula analogous to~\eqref{eq:gaugelattice}, and the measure $\mu^{\G,T}_{\YM}$ is invariant under this action. Indeed, this action preserves the reference measure $\d g$ and transforms the holonomy along loops, in this case along boundaries of faces, by conjugation, which leaves the value of the fundamental solution of the heat equation on these holonomies unchanged.
\footnote{Let us say a word about the way in which the presence of a boundary to the surface $M$ should be taken into account in \eqref{eq:discreteYM}, and how to treat the case where $M$ is not orientable. The only place where we used the orientability and orientation of $M$ is when we defined the boundary of a face as a loop winding {\em positively} around $M$. However, since the heat kernel also enjoys the invariance property $p_{t}(x)=p_{t}(x^{-1})$, it does not matter which orientation we choose around each face of the graph. Thus, \eqref{eq:discreteYM} is valid without any modification on a non-orientable surface.

In the case where $M$ has a boundary, this boundary is a finite union of circles. Our assumption that each face of a graph is homeomorphic to a disk implies that each of these circles is a path in any graph on $M$. In this case, \eqref{eq:discreteYM} still makes sense and corresponds to free boundary conditions along the boundary of $M$. Fixed boundary conditions can be imposed: it is possible to insist that the holonomy along each boundary component belongs to a specific conjugacy class in $G$. If we wish to set the boundary condition for which the holonomy along a boundary component $c=e_{1}\ldots e_{n}$ belongs to a conjugacy class $C$ of $G$, the basic ingredient is the unique probability measure $\nu_{n,C}$ on $\mathcal O_{n,C}=\{(x_{1},\ldots,x_{n})\in G^{n} : x_{n}\ldots x_{1}\in C\}$ invariant under the transitive action of $G^{n}$ given by 
\[(y_{1},\ldots,y_{n})\cdot (x_{1},\ldots,x_{n})=(y_{1}x_{1}y_{n}^{-1},y_{2}x_{2}y_{1}^{-1},\ldots,y_{n}x_{n}y_{n-1}^{-1}).\]
This measure is easily described by the formula
\[\int_{\mathcal O_{n,C}} f \d\nu_{n,C}=\int_{G^{n}} f(x_{1},\ldots,x_{n-1},x_{n}zx_{n}^{-1}x_{1}^{-1}\ldots x_{n-1}^{-1})\d x_{1}\ldots \dw x_{n},\]
for an arbitrary $z\in C$. The way in which \eqref{eq:discreteYM} should be modified is that the uniform measure on $\Conf_{\G}$ should be replaced, for the edges lying on the boundary of $M$, by the appropriate copy of a measure of the form $\nu_{n,C}$.}

\subsubsection{Invariance under subdivision}

Starting from a graph $\G$ on our surface $M$, we built the configuration space $\Conf_{\G}$ and endowed, thanks to the Driver--Sengupta formula, this space with a probability measure, the lattice $2$-dimensional Yang--Mills measure on $\G$. In doing so, we automatically produced a collection 
\[(h_{c})_{c\in \Path(\G)} \  \text{ or } \ (h_{\ell})_{\ell\in \Loop_{m}(\G)}\]
of $G$-valued random variables.\footnote{Thanks to the multiplicativity of the holonomy and the gauge invariance of the construction of the lattice Yang--Mills measure, the point of view of a collection of random variables indexed by all paths in $\G$ or by the set of loops based at a specific reference point are equivalent, see Proposition \ref{prop:multmult}.}

The property of this construction that makes it so extremely pleasant is the fact that it is {\em invariant under subdivision}.

To articulate this fundamental property, let us say that a graph $\G_{2}$ is {\em finer} than a graph $\G_{1}$ if~$\G_{2}$ can be obtained from $\G_{1}$ by subdividing and adding edges. More precisely, $\G_{2}$ is finer than $\G_{1}$ if $\E_{1}\subset \Path(\G_{2})$: each edge of $\G_{1}$ is a path in $\G_{2}$. When this happens, there is a natural map 
\begin{align*}
\Conf_{\G_{2}} & \longrightarrow \Conf_{\G_{1}} \\
g^{(2)} & \longmapsto \big(h^{(2)}_{e}(g^{(2)})\big)_{e\in \E_{1}}
\end{align*}
where each edge $e$ of $\G_{1}$ is seen as a path in $\G_{2}$ and thus assigned a holonomy by the configuration $g^{(2)}$. 

The main result of $2$-dimensional lattice Yang--Mills theory is the following.

\begin{theorem} Let $\G_{1}$ and $\G_{2}$ be two graphs on $M$. Assume that $\G_{2}$ is finer than $\G_{1}$. Then for all $T>0$, the equality $Z(\G_{1},T)=Z(\G_{2},T)$ holds and the push-forward of the measure $\mu^{\G_{2},T}_{\YM}$ by the natural map $\Conf_{\G_{2}}\to \Conf_{\G_{1}}$ is the measure $\mu^{\G_{1},T}_{\YM}$.
\end{theorem}

This theorem is so important that we are going to give an idea of the mechanism of its proof.

\begin{proof} The first observation is that one can always go from a graph to a finer graph by an appropriate succession of elementary operations consisting either in adding a new vertex in the middle of an existing edge or in adding a new edge between two existing vertices. We need to understand why neither of these elementary operations affect the partition function, nor transform essentially the measure. 

The subdivision of an edge $e$ into two new edges $e'$ and $e''$ amounts, in the integral defining the partition function and in the expression defining the discrete Yang--Mills measure, to the replacement of every occurrence of the integration variable $g_{e}$ by the product of the two new variables $g_{e''}g_{e'}$. The invariance by translation of the Haar measure ensures that this does not affect the result of any computation.

The case of the addition of a new edge is more interesting. This edge $e$ splits a face $F$ into two faces $F_{1}$ and $F_{2}$, the boundaries of which are of the form $ea$ and $be^{-1}$ for some paths $a$ and~$b$. Observe that $ba$ is a loop going along the boundary of $F$. In the computation of the partition function of the Yang--Mills measure on the finer graph, or of the integral of any functional on the configuration space of the coarser graph with respect to the image of the discrete Yang--Mills measure on the finer graph, we find an integral of a product of many factors, among which the two factors
\[p_{T|F_{1}|}\big(h_{a}(g)g_{e}\big)\ p_{T|F_{2}|}\big(g_{e}^{-1}h_{b}(g)\big)\]
contain the only two occurrences of the integration variable $g_{e}$. We can thus easily integrate with respect to $g_{e}$, using the convolution property of the heat kernel, namely the equality $p_{t}*p_{s}=p_{t+s}$, to find these two factors replaced by
\[p_{T(|F_{1}|+|F_{2}|)}\big(h_{a}(g)h_{b}(g)\big)=p_{T|F|}\big(h_{ba}(g)\big)=p_{T|F|}\big(h_{\partial F}(g)\big)\]
We are thus left with the partition function, or the integral of our functional, relative to the coarser graph.
\end{proof}
\bigskip

The partition function $Z(\G,T)$, which is now promoted to a function of $T$ alone, is a very interesting object. Let us give without proof an expression of this function. We use the notation $[a,b]=aba^{-1}b^{-1}$ for the commutator of two elements $a$ and $b$ of $G$.

\begin{proposition}\label{prop:Z} Assume that $M$ is a surface of genus $g$ without boundary. Then for all $T>0$, the partition function of the $2$-dimensional Yang--Mills theory on $M$ is given by
\[Z_{M}(T)=\int_{G^{2g}}p_{T|M|}([a_{1},b_{1}]\ldots [a_{g},b_{g}])\d a_{1}\dw b_{1}\ldots \dw a_{g} \dw b_{g}.\]
\end{proposition}

\subsubsection{The continuum limit}\label{sec:contlim}

Up to some conceptually inessential but technically annoying complications, the invariance by subdivision of the discrete theory allows one to take the limit of the discrete measures as the graphs on the surface become infinitely fine. The technical complications have to do with the fact that, because two edges of two distinct graphs can intersect in a rather pathological way, it is not always true that given two graphs, there exists a third graph that is finer than these two graphs. The net effect of this complication is the persistence, in the theorem asserting the existence and uniqueness of the Yang--Mills holonomy process, of a continuity condition. We say that a sequence of paths $(c_{n})_{n\geq 1}$ on $M$ {\em converges} to a path $c$ {\em with fixed endpoints} if all paths $c,c_{1},c_{2},\ldots$ start at the same point and finish at the same (possibly different) point, and if the sequence of the paths $(c_{n})_{n\geq 1}$ parametrised at unit speed converges uniformly to $c$.

\begin{theorem}[The Yang--Mills holonomy process, \cite{SenguptaAMS,LevyMHF}] Let $M$ be a compact surface endowed with a smooth\footnote{By a smooth measure, we mean a measure that admits a smooth positive density with respect to the Lebesgue measure in any coordinate chart.} measure of area. Let $G$ be a compact Lie group, the Lie algebra of which is endowed with an invariant scalar product. There exists a collection of $G$-valued random variables $(H_{c})_{c\in \Path(M)}$ such that
\begin{itemize}[itemsep=0pt]
\item for every graph $\G=(\V,\E,\F)$, the distribution of $(H_{e})_{e\in \E}$ is the measure $\mu^{\G,T}_{\YM}$,
\item whenever a sequence $(c_{n})_{n\geq 1}$ of paths converges with fixed endpoints to a path $c$, the sequence of random variables $(H_{c_{n}})_{n\geq 1}$ converges in probability to $H_{c}$.
\end{itemize}
Moreover, any two collections of $G$-valued random variables with these properties have the same distribution.
\end{theorem}

The Yang--Mills holonomy process $(H_{c})_{c\in \Path(M)}$ is invariant in distribution under the action of the gauge group. This means that for every function $g:M\to G$, the following equality in distribution holds:
\begin{equation}\label{eq:GI}
\Big(g(\overline{c})^{-1}H_{c}g(\underline{c})\Big)_{c\in \Path(M)}\stackrel{(d)}{=}(H_{c})_{c\in \Path(M)}
\end{equation}
where $\underline{c}$ and $\overline{c}$ denote respectively the starting and finishing point of a path $c$.
In particular, the distribution of $H_{c}$ is uniform on $G$ for every path $c$ that is not a loop. Of course, this huge collection of uniform random variables is correlated in a complicated way, in particular to allow the random variables associated with loops to have non-uniform distributions.

The holonomy process also enjoys a property of invariance under area-preserving maps of $M$: if $\phi:M\to M$ is an area-preserving diffeomorphism, then $\phi$ preserves the class $\Path(M)$ and the family $(H_{\phi(c)})_{c\in \Path(M)}$ has the same distribution as the family $(H_{c})_{c\in \Path(M)}$. This is because the Driver--Sengupta formula depends only on the combinatorial structure of the graph under consideration, and on the areas of its faces. This is consistent with the fact that the Yang--Mills action, which we originally defined on a Riemannian manifold by \eqref{eq:defSint}, depends, if the manifold is $2$-dimensional, on the Riemannian structure only through the Riemannian area. We already mentioned this important point in relation with the expression \eqref{eq:YMH} of the Yang--Mills action. 

\subsubsection{The structure of the holonomy process}\label{sec:structYMHP}

The structure of the Yang--Mills holonomy process can be described fairly concretely provided one understands the structure of the set of loops on a graph. 

Let us consider a graph $\G$ on $M$ and a vertex $m$ of this graph. We denote naturally by $\Loop_{m}(\G)$ the set of loops in $\G$ based at $m$. The operation of concatenation makes $\Loop_{m}(\G)$ a monoid, with unit element the constant loop at $m$. Each element $\ell$ of this monoid has an `inverse' $\ell^{-1}$, but it is not true, unless $\ell$ is already the constant loop, that $\ell\ell^{-1}$ is the constant loop. In order to make $\Loop_{m}(\G)$ a group, into which $\ell^{-1}$ is truly the inverse of $\ell$, it is natural to introduce on it the {\em backtracking equivalence} relation, for which two loops are equivalent if one can go from one to the other by successively erasing or inserting sub-loops of the form $ee^{-1}$, where $e$ is an edge of the graph. 

Each equivalence class of loops contains a unique loop of shortest length, which is also the unique {\em reduced} loop in this class, where by a reduced loop we mean one without any sub-loop of the form $ee^{-1}$.

Moreover, concatenation is compatible with this equivalence relation and the quotient monoid is a group. This quotient monoid can be more concretely described as the set $\Loop_{m}^{\rm red}(\G)$ of reduced loops endowed with the operation of concatenation-followed-by-reduction. 

With this group of reduced loops in hand, we can make several observations. 

\begin{itemize}[itemsep=0pt] 
\item Each element $g$ of the configuration space $\Conf_{\G}$ induces, by the holonomy map, a map $\Loop_{m}^{\rm red}(\G)\to G$, which sends a loop $\ell$ to $h_{\ell}(g)$. This map is a group homomorphism, and the map
\[\Conf_{\G} \longrightarrow {\rm Hom}(\Loop_{m}^{\rm red}(\G),G)\]
is onto. Moreover, it descends to a bijection 
\[\Conf_{\G}/{\rm Maps}(\V,G) \build{\longrightarrow}_{}^{\sim} {\rm Hom}(\Loop_{m}^{\rm red}(\G),G)/G\]
where the action on the left is that of the gauge group, and on the action on the right is that of $G$ by conjugation.
\item Let $\Gamma$ denote the $1$-skeleton of the graph, that is, the union of the ranges of its edges. The map $\Loop_{m}^{\rm red}(\G)\to \pi_{1}(\Gamma,m)$ which simply sends a reduced loop to its homotopy class is an isomorphism. 
\item The group $\Loop_{m}^{\rm red}(\G)$, being isomorphic to the fundamental group of a graph, or of a $1$-dimensional complex, is a free group. The rank of this group is equal to $|\E|-|\V|+1=|\F|-\chi(M)+1=|\F|+2g-1$, where $\chi(M)$ is the Euler characteristic of $M$ and $g$ its genus.
\end{itemize}

It is useful to recognise that the free group $\Loop_{m}^{\rm red}(\G)$ admits nice bases.\footnote{Recall that a free group admits bases, that is, subsets by which it is freely generated. Any two bases have the same cardinality, called the rank of the group. Any subgroup of a free group is itself a free group, but the rank of a subgroup can be larger than the rank of the group. In fact, the free group of rank $2$ contains subgroups of arbitrary finite or (countably) infinite rank.} Let us call {\em lasso} around a face $F$ of $\G$ any loop of the form $c.\partial F .c^{-1}$, where $c$ is a path from $m$ to a vertex on the boundary of $F$, and $\partial F$ is a loop going once around $F$. 

It is now quite easy to describe the holonomy process. Let us begin with the case of the plane, or the disk.

\begin{proposition} Assume that $M$ is a disk or the plane. Let $\G$ be a graph on $M$. The free group $\Loop_{m}^{\rm red}(\G)$ admits a basis $\{\lambda_{F}: F\in \F\}$ such that
\begin{itemize}[itemsep=0pt]
\item for each face $F$, the loop $\lambda_{F}$ is a lasso around $F$,
\item under the lattice Yang--Mills measure $\mu^{\G,T}_{\YM}$, the random variables $(H_{\lambda_{F}} : F\in \F)$  are independent, each $H_{\lambda_{F}}$ being distributed according to the measure $p_{T|F|}(g)\d g$.
\end{itemize}
\end{proposition}

In a sense, the holonomy process has independent increments distributed according to the fundamental solution of the heat equation: it can be described as a `Brownian motion on $G$ indexed by loops' on the disk, or on the plane. The role of time is played by area, and increments occur along faces of the graph, or lassos, instead of intervals of time.

In the case of a closed surface, the situation is slightly different. In this case, the most natural presentation of the group $\Loop_{m}^{\rm red}(\G)$ is not as a free group (which it is), but with one generator too many, and one relation. 

\begin{proposition} Assume that $M$ is a closed surface of genus $g$. Let $\G$ be a graph on $M$. Set $r=|\F|$. The free group $\Loop_{m}^{\rm red}(\G)$ admits a presentation
\[\Loop_{m}^{\rm red}(\G)=\big\langle \lambda_{F_{1}},\ldots,\lambda_{F_{r}}, a_{1},b_{1},\ldots,a_{g},g_{b} \ \big|\  [a_{1},b_{1}]\ldots [a_{g},b_{g}]=\lambda_{F_{1}}\ldots \lambda_{F_{r}}\big\rangle\]
where
\begin{itemize}[itemsep=0pt]
\item the loops $\lambda_{F_{1}},\ldots,\lambda_{F_{r}}$ are lassos around the $r$ faces of $\G$,
\item the homotopy classes of the loops $a_{1},b_{1},\ldots,a_{g},b_{g}$ generate  $\pi_{1}(M,m)$,
\item for every test function $f:G^{2g+r}\to \C$, one has
\begin{align}\label{eq:distYM}
\int_{\Conf_{\G}} f(H_{\lambda_{1}},\ldots,H_{\lambda_{r}},H_{a_{1}},H_{b_{1}},\ldots,H_{a_{g}},H_{b_{g}}) \d \mu^{\G,T}_{\YM}&\\ \nonumber
&\hspace{-7.3cm} =Z_{M}(T)^{-1}\int_{G^{2g+r-1}}f(z_{1},\ldots,z_{r-1},z_{r},x_{1},y_{1},\ldots,x_{g},y_{g}) p_{T|F_{1}|}(z_{1})\ldots p_{T|F_{r}|}(z_{r})\\ \nonumber
&\hspace{0.45cm} \d z_{1}\ldots \d z_{r-1} \d x_{1} \d y_{1} \ldots \d x_{g} \d y_{g}
\end{align}
where in the last integral, $z_{r}$ stands for
\[z_{r}=(z_{r-1}\ldots z_{1}[a_{g},b_{g}]\ldots [a_{1},b_{1}])^{-1}\]
\end{itemize}
\end{proposition}

Let us try to spell out the probabilistic content of this result. The presentation of the group $\Loop_{m}^{\rm red}(\G)$ that we chose splits it into a homotopically trivial part, giving rise to the random variables $H_{\lambda_{1}},\ldots,H_{\lambda_{r}}$, and a system of generators of the fundamental group of $M$, associated with the random variables $H_{a_{1}},H_{b_{1}},\ldots,H_{a_{g}},H_{b_{g}}$. A particular role is played by the homotopically trivial loop $C=[a_{1},b_{1}]\ldots [a_{g},b_{g}]$.
\begin{itemize}[itemsep=0pt]
\item The distribution of the random variable $H_{C}$ is such that for every continuous test function $\tilde f:G\to \C$, 
\[\int_{\Conf_{G}} \tilde f(H_{C})\d\mu^{\G,T}_{\YM}=Z_{M}(T)^{-1}\int_{G^{2g}}(\tilde fp_{T|M|})([a_{1},b_{1}]\ldots [a_{g},b_{g}])\d a_{1} \d b_{1}\ldots \d a_{g} \d b_{g}\]
This does not seem to be a particularly well-known distribution. It needs not have a density with respect to the Haar measure: for instance if $G=\U(N)$, it is supported by the Haar-negligible subgroup $\SU(N)$. However, it is, by definition, absolutely continuous with respect to the distribution of the product of $g$ independent commutators of independent uniformly distributed random variables, and this distribution, for example if $G=\SU(N)$ and provided $g\geq 2$, is absolutely continuous with respect to the Haar measure. It is also possible to write a Fourier series for this distribution, but it involves Littlewood--Richardson coefficients, or more generally an understanding of the tensor product of irreducible representations of $G$.
\item Conditional on $H_{C}$, the families $(H_{\lambda_{1}},\ldots,H_{\lambda_{r}})$ and $(H_{a_{1}},H_{b_{1}},\ldots,H_{a_{g}},H_{b_{g}})$ are independent. It is also true that the random variables
\[(H_{\lambda_{1}},\ldots,H_{\lambda_{r}})\text{ mod }G \ \ \ \text{ and } \ \ \ (H_{a_{1}},H_{b_{1}},\ldots,H_{a_{g}},H_{b_{g}})\text{ mod }G\]
with values in $G^{r}/G$ and $G^{2g}/G$, where $G$ acts by conjugation, are independent conditional on $H_{C}\text{ mod }G$, that is, conditional on the conjugacy class of $H_{C}$. 
\end{itemize}

On a surface of genus $g$, the probabilistic backbone of the holonomy process can thus be described as consisting of a segment of a Brownian motion on $G$ of length $T|M|$ and $2g$ independent Haar distributed random variables on $G$, jointly conditioned on the final point of the Brownian motion being equal to the products of the $g$ commutators of the uniform random variables taken in pairs.

The case where $M$ is a sphere is special, in the sense that it involves no uniform random variables, but only a Brownian bridge on $G$ going from $1_{G}$ to $1_{G}$ in a time equal to $T$ times the total area of the sphere.

\subsection{Wilson loop expectations}

A different approach to the description of the distribution of the Yang--Mills holonomy process consists in identifying a natural class of scalar, gauge invariant, functionals of this process, the distribution of which is hoped to contain as much information as possible. The most natural class of such functionals is that of {\em Wilson loop functionals}, which are indeed the most important scalar observables of the theory. A Wilson loop functional is constructed by choosing a certain number of loops $\ell_{1},\ldots,\ell_{n}$  on $M$, then the same number of conjugation-invariant functions $\chi_{1},\ldots,\chi_{n}:G\to \C$ and by forming the product
\begin{equation}\label{eq:WLF}
\chi_{1}(H_{\ell_{1}})\ldots \chi_{n}(H_{\ell_{n}})
\end{equation}
When $G$ is a group of matrices, the simplest choice of conjugation-invariant function is the trace. The {\em Wilson loop expectations}, which play in this theory the role of $n$-point functions, are the numbers
\begin{equation}\label{eq:WLE}
\E[\Tr(H_{\ell_{1}})\ldots \Tr(H_{\ell_{n}})]
\end{equation}
the computation of which is a seemingly endless subject of reflection. We will discuss in the next section a few concrete examples of computation of such numbers. For the time being, let us say a word about the amount of information that they carry. 

Suppose we know the collection of all the numbers \eqref{eq:WLE}, or more generally the expectation of all functionals of the form \eqref{eq:WLF}. Then we know the joint distribution of all random variables of the form $\chi(H_{\ell})$ where $\ell$ is a loop and $\chi:G\to \C$ is an invariant function. Since $G$ is compact, invariant functions separate conjugacy classes and we know, in fact, the joint distribution of the conjugacy classes of all variables $H_{\ell}$. This is certainly an important piece of information. However, the form of the action of the group of gauge transformations on the collection of holonomies, as given by \eqref{eq:GI}, indicates that this action preserves more than just the individual conjugacy classes of the holonomies. Indeed, if $\ell_{1},\ldots,\ell_{n}$ are based at the same point, then it is the orbit of $(H_{\ell_{1}},\ldots,H_{\ell_{n}})$ under the operation of {\em simultaneous conjugation}
\[(h_{1},\ldots,h_{n})\mapsto (gh_{1}g^{-1},\ldots,gh_{n}g^{-1})\]
that is gauge-invariant. To grasp the geometric meaning of this invariance, it is useful to take a concrete example for $G$, say $G=\SU(N)$ or even $G=\SO(3)$. In these groups, knowing the individual conjugacy classes of a collection of elements amounts to knowing their eigenvalues, that is, in the case of $\SO(3)$, the angles of the rotations. On the other hand, to know the orbit of these elements under simultaneous conjugation requires the additional knowledge of the relative positions of their eigenspaces, or for rotations, the relative positions of their axes.

The main question is then the following. Is it the case that the Wilson loop expectations describe not only the individual conjugacy classes of the $G$-valued random variables that constitute the Yang--Mills process, but also the simultaneous conjugacy class of all variables associated to the loops based at some point $m$ of $M$? In more precise terms, is it true that the algebra of functions on $\A/\J$ generated by Wilson loop functionals separates points? If not, it cannot be said that the Wilson loop functionals constitute a complete set of gauge-invariant scalar observables. 

The answer turns out to depend entirely on the group $G$, and it does not seem to be known in all cases, even for compact Lie groups.\footnote{It would be more prudent to say that it is not known to the author.} The property that $G$ must have for the answer to be positive is the following.\footnote{The name of Property W is by no means standard.}

\begin{definition}[Property W] We say that a group $G$ has the property W if for any $n\geq 2$ and any two collections $x_{1},\ldots,x_{n}$ and $x'_{1},\ldots,x'_{n}$ of elements of $G$, the assumption that every word in $x_{1},\ldots,x_{n}$ and their inverses  is conjugated to the same word in $x'_{1},\ldots,x'_{n}$ and their inverses implies the existence of an element $y$ of $G$ such that $x'_{1}=yx_{1}y^{-1}, \ldots, x'_{n}=yx_{n}y^{-1}$.
\end{definition}

Since this long definition is maybe not very pleasant to read, let us word it differently. We are comparing two relations between $n$-tuples $(x_{1},\ldots,x_{n})$ and $(x'_{1},\ldots,x'_{n})$ of elements of $G$. The first is the relation of simultaneous conjugation
\begin{equation}\label{eq:SC} \tag{SC}
\exists y\in G, \ x'_{1}=yx_{1}y^{-1}, \ldots, x'_{n}=yx_{n}y^{-1}
\end{equation}
The second could be called {\em lexical conjugation} and holds exactly when
\begin{equation}\label{eq:LC} \tag{LC}
\text{every word in $x_{1},\ldots,x_{n}$ is conjugated to the same word in $x'_{1},\ldots,x'_{n}$}
\end{equation}
where a word in a certain set of letters can involve these letters and their inverses. We also considered a third property of individual conjugation
\begin{equation}\label{eq:IC} \tag{IC}
\exists y_{1},\ldots,y_{n} \in G, \ x'_{1}=y_{1}x_{1}y_{1}^{-1}, \ldots, x'_{n}=y_{n}x_{n}y_{n}^{-1}
\end{equation}
In any group, one has the chain of implications
\[\eqref{eq:SC} \Rightarrow \eqref{eq:LC} \Rightarrow \eqref{eq:IC}\] 
Unless the group $G$ has very special properties (for instance that of being abelian), the second implication is not an equivalence, and the property \eqref{eq:IC} is much weaker than the property \eqref{eq:LC}. For the group~$G$ to have the property W means that the properties \eqref{eq:SC} and \eqref{eq:LC} are equivalent. The proof of the following result can be found in \cite{LevySN}, see also \cite{Durhuus,SenguptaGI}.

\begin{theorem} Any Cartesian product of special orthogonal, orthogonal, special unitary, unitary, and symplectic groups has the property W.
\end{theorem}

It is known that some non-compact groups fail to have the property W. However, it seems not be known wether this equivalence holds, for instance, for spin groups. 

\section{Computation of Wilson loop expectations}\label{sec:compWLE}

In this section, we will give a few concrete examples of computations with the Yang--Mills holonomy process, with an eye to its so-called {\em large $N$ limit}, that is, its behaviour when the group $G$ is taken to be $\U(N)$ with an appropriately scaled invariant product on its Lie algebra, and $N$ tends to infinity.\footnote{The notion of large $N$ limit also applies to the cases where $G=\SO(N)$ and $G=\Sp(N)$, the real and quaternionic analogues of $\U(N)$ or $\SU(N)$. As far as we understand today, there is no essential difference between the three cases. More precisely, the computations for finite $N$ are similar in the three cases, if generally a bit more complicated in the orthogonal case and even more so in the symplectic case, and the large $N$ limits are identical. }

The basis of virtually any computation in $2$-dimensional Yang--Mills theory is the Driver--Sengupta formula \eqref{eq:discreteYM}. This formula can be combined with an expression of the heat kernel on~$G$, for example its Fourier expansion, and lead to very concrete calculations.  It is also possible to use a more dynamical, either analytic or probabilistic approach to the heat kernel, by seeing it as the solution of the heat equation or, almost equivalently, as the density of the distribution of the Brownian motion on $G$. We will illustrate these possibilities on a few examples in the simplest case where $M$ is the plane, and then turn to the much more complicated case where $M$ is the $2$-dimensional sphere. For the sake of simplicity, we will assume in this section that the coupling constant $T$ that appears in \eqref{eq:discreteYM}  is equal to $1$.

\subsection{The Brownian motion on the unitary group}

In order to be as concrete as possible, and because we are interested in the large $N$ limit, we will in this section choose $G=\U(N)$, the unitary group of rank $N$. As indicated earlier (see Footnote \ref{fn:GU}), we endow the Lie algebra of $\U(N)$, which is the space $\u(N)$ of $N\times N$ skew-Hermitian matrices, with the scalar product $\langle X,Y\rangle=N\Tr(X^{*}Y)$. In the Euclidian space $(\u(N),\langle \cdot, \cdot \rangle)$, we consider a linear Brownian motion $(K_{t})_{t\geq 0}$, use it to form the stochastic differential equation
\begin{equation}\label{eq:EDS}
dU_{t}=U_{t}\d K_{t} - \frac{1}{2}U_{t}\d t\ , \ \ U_{0}=I_{N}
\end{equation}
and call the unique solution to this equation the Brownian motion on $\U(N)$. 

Using the notation $\Tr$ for the usual trace of a $N\times N$ matrix and $\tr=\frac{1}{N}\Tr$ for its normalised trace, the usual rules of stochastic calculus take, in this matricial context, the following nice form: for all $N\times N$ matrix $A$, measurable with respect to $\sigma(K_{s} : s\leq t)$, we have
\begin{equation}\label{eq:KK}
\d K_{t} A \d K_{t}  = - \tr(A) \d t \ \ \text{ and } \ \  \d K_{t} \, \tr(A \d K_{t}) =-\frac{1}{N^{2}} A \d t
\end{equation}
This relation can be used to check that $\d(U_{t}U_{t}^{*})=0$, so that the trajectories of the process $B$ stay almost surely, as expected, in $\U(N)$.

The density of the distribution of $U_{t}$ with respect to the normalised Haar measure on $\U(N)$ is the function $p_{t}$ appearing in the Driver--Sengupta formula, and that we described in Section~\ref{sec:DS}.

It will be useful to know the Fourier series of this function $p_{t}:\U(N)\to \R$. To describe it, let us introduce the set $\widehat\U(N)$ of equivalence classes of irreducible representations (or {\em irreps}) of~$\U(N)$. For every $\alpha\in \widehat \U(N)$, let us denote by $d_{\alpha}$ the degree of $\alpha$, that is, the dimension of the space on which $\U(N)$ acts through $\alpha$. Let us also denote by $\chi_{\alpha}:\U(N)\to \C$ the character of~$\alpha$, and by $c_{2}(\alpha)$ the quadratic Casimir number of $\alpha$, that is, the non-negative real number such that 
\[\Delta \chi_{\alpha}=-c_{2}(\alpha)\chi_{\alpha}\]
The Fourier series of the heat kernel is then
\begin{equation}\label{eq:Fourierpt}
p_{t}=\sum_{\alpha\in \widehat \U(N)} e^{-\frac{c_{2}(\alpha)t}{2}} d_{\alpha}\chi_{\alpha}
\end{equation}
and there is nothing specific to $\U(N)$ in this formula. 

It is however possible, in the case of $\U(N)$, to write explicitly each of its ingredients. Indeed, the set of irreps of $\U(N)$ is conveniently labelled by non-increasing sequences of $N$ relative integers $\lambda=(\lambda_{1}\geq \ldots \geq \lambda_{N})$, called dominant weights. The dimension and quadratic Casimir number of the irrep with highest weight $\lambda$ are given by the formulas
\begin{equation}\label{eq:dimcas}
d_{\lambda}=\prod_{1\leq i< j \leq N} \frac{\lambda_{i}-\lambda_{j}+j-i}{j-i} \ \text{ and } \ N c_{2}(\lambda)=\sum_{1\leq i \leq N} \lambda_{i}^{2}+\sum_{1\leq i < j \leq N} (\lambda_{i}-\lambda_{j})
\end{equation}
The character of this representation is given, up to a power of the determinant, by a Schur function, but we will not need its explicit formula. 

We are now equipped to make some computations with the Yang--Mills holonomy process.

\subsection{The simple loop on the plane} \label{sec:simpleplane}

\subsubsection{Using harmonic analysis}

\begin{figure}[h!]
\begin{center}
\rotatebox{45}{\includegraphics[width=4cm]{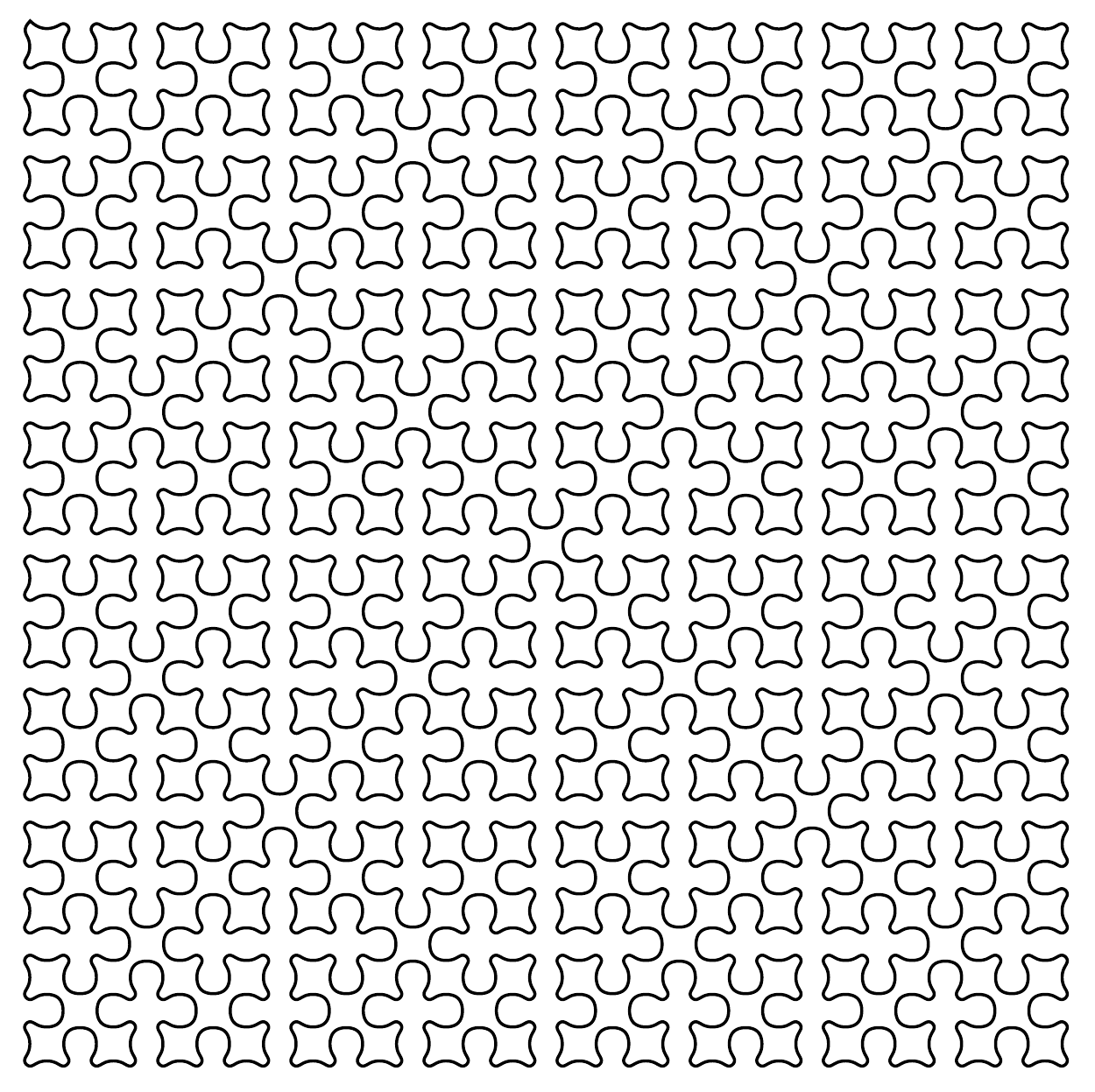}}
\caption{\label{fig:simple} \small A simple loop on the plane}
\end{center}
\end{figure}

Let us consider, on the plane, a loop $\ell$ that is a simple loop going once around a domain of area $t$ (see, if needed, Fig. \ref{fig:simple}). The partition function of the Yang--Mills model on the plane is equal to $1$ and the Driver--Sengupta formula \eqref{eq:discreteYM} tells us that for every continuous test function $f:\U(N)\to \C$, we have
\[\E[f(H_{\ell})]=\int_{\U(N)} f(x)p_{t}(x)\d x \]
In other words, $H_{\ell}$ has the same distribution as $U_{t}$, the value at time $t$ of the Brownian motion on $\U(N)$ defined in the previous section.

Using the Fourier expansion \eqref{eq:Fourierpt} and the classical orthogonality relations between characters, we find, for every irrep $\alpha$ of $\U(N)$ acting on the vector space $V_{\alpha}$, the equality
\[\E[\alpha(H_{\ell})]=e^{-\frac{c_{2}(\alpha)t}{2}}\, \id_{V_{\alpha}}\]
which holds in ${\rm End}(V_{\alpha})$.
In particular, since the usual trace is, on $\U(N)$, the character of the natural representation, which has highest weight $(1,0,\ldots,0)$, dimension $N$ and quadratic Casimir~$1$, we find
\begin{empheq}[box=\mybox]{equation}\label{eq:tr}
\E[H_{\ell}]=e^{-\frac{t}{2}} I_{N} \ \ \text{ and } \ \ \E[\tr (H_{\ell})]=e^{-\frac{t}{2}}
\end{empheq}

Suppose now that we want to compute the expectation of $\tr(H_{\ell}^{2})$, which is also the expectation of $\tr(H_{\ell^{2}})$, where $\ell^{2}$ is the loop $\ell$ gone along twice. From the Driver--Sengupta formula and the Fourier expansion of the heat kernel, we get the expression
\[\E[\tr(H_{\ell}^{2})]=\sum_{\lambda\in \widehat\U(N)} e^{-\frac{c_{2}(\lambda)t}{2}} d_{\lambda} \int_{\U(N)} \tr(x^{2})\chi_{\lambda}(x)\d x\]
In order to go further, we need to know that, at least when $N\geq 2$,
\[\tr(x^{2})=\chi_{(2,0,\ldots,0)}(x)-\chi_{(1,1,0\ldots,0)}(x)\]
Using again the orthogonality of characters, we find, after some reordering of the terms,
\begin{empheq}[box=\mybox]{equation} \label{eq:tr2}
\E[\tr(H_{\ell}^{2})]=e^{-t} \Big(\cosh \frac{t}{N}-N\sinh \frac{t}{N}\Big)
\end{empheq}

It is possible to go further down this road, by systematically writing the function $x\mapsto \tr(x^{n})$ as a linear combination of characters. This is what Philippe Biane did to determine the large $N$ limit of the non-commutative distribution of the Brownian motion on the unitary group. The simplest non-trivial case is the large $N$ limit of \eqref{eq:tr2}:
\begin{empheq}[box=\mybox]{equation} \label{eq:tr2N}
\lim_{N\to \infty}\E[\tr(H_{\ell}^{2})]=e^{-t}(1-t)
\end{empheq}

The general formula is nice enough, at least in the limit when $N$ tends to infinity, to be quoted explicitly. It was discovered independently by Philippe Biane and Eric Rains, who formulated it in terms of the Brownian motion on $\U(N)$ rather than the Yang--Mills holonomy process.

\begin{theorem}[Biane \cite{Biane}, Rains \cite{Rains}] With the current notation, and for every integer $n\geq 1$, 
\begin{equation}\label{eq:limmomMB}
\lim_{N\to \infty} \E[\tr(H_{\ell}^{n})]=e^{-\frac{nt}{2}}\sum_{k=0}^{n-1} \frac{(-t)^{k}}{k!} n^{k-1}\binom{n}{k+1}
\end{equation}
\end{theorem}

It must be said that this result already appeared, without proof, in Isadore Singer's seminal paper on the large $N$ limit of the Yang--Mills holonomy field \cite{Singer}.\footnote{Singer and Rains recognise, in the right-hand side of \eqref{eq:limmomMB}, modified Laguerre polynomials of the first kind. As far as I know, a structural explanation for the appearance of these polynomials in this context has yet to be given.}
 
One of Biane's aims in \cite{Biane} was to prove the following theorem concerning the limit as~$N$ tends to infinity of the Brownian motion on $\U(N)$ as a stochastic process. This convergence result is stated in the language of free probability, a theory presented in detail in the book of Alexandru Nica and Roland Speicher \cite{NicaSpeicher}.

\begin{theorem}[Biane \cite{Biane}] As $N$ tends to infinity, the Brownian motion on $\U(N)$ converges in non-commutative distribution, as a process,  towards a unitary non-commutative process $(u_{t})_{t\geq 0}$ with free stationary multiplicative increments such that for all integer $n\geq 0$ and all real $t\geq 0$, the expectation of~$u_{t}^{n}$ and that of $(u_{t}^{*})^{n}$ are given by the right-hand side of \eqref{eq:limmomMB}.

\end{theorem}

\subsubsection{Using stochastic calculus}

Let us illustrate, on the same example of a simple loop on the plane, the dynamical approach to the same computations, based on the use of It\={o}'s formula. The general principle of these computations is to see the quantities such as the left-hand sides of \eqref{eq:tr} and \eqref{eq:tr2} as functions of $t$, and to write a differential equation that they satisfy. Recall that $t$, in our current notation, is the area of the disk enclosed by the simple loop $\ell$. A variation of $t$ can thus be described, in geometrical terms, as a variation of the area of the unique face enclosed by $\ell$.

As a first example, let us use \eqref{eq:EDS} and It\={o}'s formula to find
\[\frac{d}{dt}\E[\tr(H_{\ell})]=\frac{d}{dt}\E[\tr(U_{t})]=-\frac{1}{2}\E[\tr(U_{t})]\]
which, together with the information $\E[\tr(U_{0})]=1$, yield immediately \eqref{eq:tr}.

Let us apply the same strategy to the computation of $\E[\tr(H_{\ell}^{2})]=\E[\tr(U_{t}^{2})]$. The computation is more interesting and involves the first of the two rules \eqref{eq:KK}. We find
\begin{equation}\label{eq:sys1}
\frac{d}{dt}\E[\tr(U_{t}^{2})]=-\E[\tr(U_{t}^{2})]-\E[\tr(U_{t})^{2}]
\end{equation}
and see a function of $t$ pop up that we were initially not interested in, namely $\E[\tr(U_{t})^{2}]$.
The only way out left to us is retreat forwards and we compute the derivative with respect to $t$ of this new function, using now the second rule of \eqref{eq:KK}:
\begin{equation}\label{eq:sys2}
\frac{d}{dt}\E[\tr(U_{t})^{2}]=-\frac{1}{N^{2}}\E[\tr(U_{t}^{2})]-\E[\tr(U_{t})^{2}]
\end{equation}
All's well that ends well: \eqref{eq:sys1} and \eqref{eq:sys2} form a closed system of ordinary differential equations that is easily solved and from which we recover, in particular, \eqref{eq:tr2}. As a bonus, we get
\begin{empheq}[box=\mybox]{equation} \label{eq:tr22}
\E[\tr(H_{\ell})^{2}]=e^{-t} \Big(\cosh \frac{t}{N}-\frac{1}{N}\sinh \frac{t}{N}\Big)
\end{empheq}
The only change with respect to \eqref{eq:tr2} is the change from $N$ to $\frac{1}{N}$ in front of the hyperbolic sine, with the effect that
\begin{equation}\label{eq:tr22fac}
\lim_{N\to \infty}\E[\tr(H_{\ell})^{2}]=e^{-t}=\lim_{N\to \infty}\E[\tr(H_{\ell})]^{2}
\end{equation}
This is an instance of a general factorisation property which was observed, among others, by Feng Xu \cite{Xu}, and which is a consequence of the concentration, in the limit where $N$ tends to infinity, of the spectra of the random matrices that we are considering.

\subsection{Yin~.~.~.}\label{sec:Yin}

\begin{figure}[h!]
\begin{center}
\includegraphics{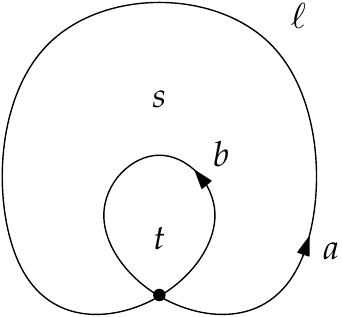}
\caption{\label{fig:heart} \small The loop $\ell$ goes first once along the larger circle (the edge $a$) and then once along the smaller circle (the edge $b$). The loop $ab$ is equivalent to the concatenation of~$ab^{-1}$, $b$ and $b$. The loops $ab^{-1}$ and $b$ are essentially simple loops surrounding disjoint domains.}
\end{center}
\end{figure}

Let us consider a slightly more complicated loop depicted on Fig. \ref{fig:heart}. This loop goes once around a domain of area $s+t$ and then once around a smaller domain of area $t$ contained in the first one. 

Let us apply the Driver--Sengupta formula in this case. We denote a generic element of the configuration space $\U(N)^{2}$ by $(x_{a},x_{b})$, in relation with our labelling by $a$ and $b$ of the two edges of the graph formed by $\ell$. Thus, for every continuous test function $f:\U(N)\to \C$, we have
\[\E[f(H_{\ell})]=\int_{\U(N)^{2}}f(x_{b}x_{a}) p_{s}(x_{b}^{-1}x_{a})p_{t}(x_{b})\d x_{a}\d x_{b}\]
Note that, according to \eqref{eq:defholo}, the discrete holonomy map is order-reversing, so that the loop $\ell=ab$ gives rise to the map $h_{\ell}(x_{a},x_{b})=x_{b}x_{a}$.

The change of variables $(y,z)=(x_{b}^{-1}x_{a},x_{b})$ preserves the Haar measure on $\U(N)^{2}$ and we have
\begin{equation}\label{eq:expheart}
\E[f(H_{\ell})]=\int_{\U(N)^{2}}f(z^{2}y) p_{s}(y)p_{t}(z)\d y \d z
\end{equation}
 This corresponds to the fact, explained in the caption of Fig. \ref{fig:heart}, that the loop $\ell$ can be written as $\ell_{1}\ell_{2}\ell_{2}$, where $\ell_{1}$ goes around the moon-shaped domain sitting between the two disks, and $\ell_{2}$ goes around the small circle of area $t$. These loops enclose disjoint domains, and although $\ell_{1}$ is not strictly speaking self-intersection free, they are essentially simple, in the sense that they can be approximated by simple loops. 

From this graphical decomposition of $\ell$, or from \eqref{eq:expheart}, we infer that $H_{\ell}$ has the distribution of~$V_{t}^{2}U_{s}$, where $U$ and $V$ are independent Brownian motions on $\U(N)$.\footnote{Thanks to the independence of the multiplicative increments of the Brownian motion, this distribution is of course also that of $V_{t}^{2}(V_{t}^{-1}V_{t+s})=V_{t}V_{s+t}$. Reasoning in this way amounts to undo the change of variables that we did to obtain \eqref{eq:expheart}.} Using the independence, the fact that the expectation of $U_{s}$ is $e^{-\frac{s}{2}}I_{N}$ (see \eqref{eq:tr}), and \eqref{eq:tr2}, we find

\begin{empheq}[box=\mybox]{equation} \label{eq:trheart}
\E[\tr(H_{\ell})]=e^{-\frac{s}{2}-t} \Big(\cosh \frac{t}{N}-N\sinh \frac{t}{N}\Big)
\end{empheq}
and, letting $N$ tend to infinity,
\begin{empheq}[box=\mybox]{equation} \label{eq:trheartMF}
\lim_{N\to\infty}\E[\tr(H_{\ell})]=e^{-\frac{s}{2}-t} (1-t)
\end{empheq}

We succeeded in computing the expectation of $\tr(H_{\ell})$, but we did so by taking advantage of the favourable circumstances, namely the fact that the word $V_{t}^{2}U_{s}$ is a very simple one, with two independent Brownian motions appearing one after the other (and not, for example, as $U_{s}V_{t}U_{s}V_{t}$), and the fact that the expectation of $U_{s}$ is a very simple matrix.

A more systematic approach is possible, by looking at $\E[\tr(V_{t}^{2}U_{s})]$ as a function of $s$
 and $t$ and by using It\={o}'s formula to compute its partial derivatives. One finds
 \begin{align*}
& \partial_{s}\E[\tr(V_{t}^{2}U_{s})]=-\frac{1}{2}\E[\tr(V_{t}^{2}U_{s})]\\
&\partial_{t}\E[\tr(V_{t}^{2}U_{s})]=-\E[\tr(V_{t}^{2}U_{s})]-\E[\tr(V_{t})\tr(V_{t}U_{s})]
\end{align*}
Once again, a function appears that we were not considering at first. Let us apply the same treatment to this new function:
 \begin{align*}
& \partial_{s}\E[\tr(V_{t})\tr(V_{t}U_{s})]=-\frac{1}{2}\E[\tr(V_{t})\tr(V_{t}U_{s})]\\
&\partial_{t}\E[\tr(V_{t})\tr(V_{t}U_{s})]=-\E[\tr(V_{t})\tr(V_{t}U_{s})]-\frac{1}{N^{2}}\E[\tr(V_{t}^{2}U_{s})]
\end{align*}
It is possible to solve this system and to recover \eqref{eq:trheart}. 

An interesting observation is the fact that the linear combination $2\partial_{s}-\partial_{t}$ of partial derivatives is particularly simple:
\begin{align} \label{eq:eyeMM1}
&(2\partial_{s}-\partial_{t})\E[\tr(V_{t}^{2}U_{s})]=\E[\tr(V_{t})\tr(V_{t}U_{s})]\\
\text{and } & (2\partial_{s}-\partial_{t})\E[\tr(V_{t})\tr(V_{t}U_{s})]=\frac{1}{N^{2}}\E[\tr(V_{t}^{2}U_{s})]
\end{align}
These are instances of the Makeenko--Migdal equations that we will discuss in greater detail in the next section. Before that, let us study another example.

\subsection{.~.~.~and Yang}

\begin{figure}[h!]
\begin{center}
\includegraphics{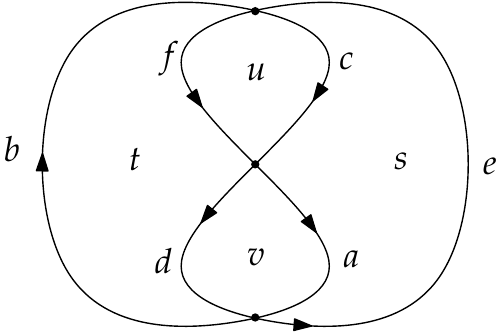}
\caption{\label{fig:eight} \small An eight-shaped loop on the plane. The letters $s,t,u,v$ in the faces indicate the areas of the faces. The other letters label the edges of the graph. The loop can be decomposed, as we did for the heart-shaped loop, as a product of lassos enclosing paiwise disjoint domains: 
$abcdef=(ad^{-1}) (dbf)(f^{-1}c)(da^{-1})(aec)(c^{-1}f)$. Here, by a lasso, we mean a loop of the form $clc^{-1}$, where $c$ is a path starting from the starting point of our loop and $l$ is a simple loop. In this particular case, the path $c$ is always the constant path.
}
\end{center}
\end{figure}

Let us now consider the eight-shaped loop drawn on Fig. \ref{fig:eight}. The Driver--Sengupta formula yields, with the by now usual notation, and taking the inversion of the order into account,
\[\E[f(H_{\ell})]=\int_{\U(N)^{6}} f(x_{f}x_{e}x_{d}x_{c}x_{b}x_{a}) p_{s}(x_{a}x_{c}x_{e})p_{t}(x_{f}x_{b}x_{d})p_{u}(x_{c}^{-1}x_{f})p_{v}(x_{a}^{-1}x_{d})\d x\]

The appropriate change of variables is dictated by the geometry of the loop, more precisely by a decomposition in product of lassos, one of which is given in the caption of Fig. \ref{fig:eight}. Accordingly, let us set
\[(y,z,g,h,e,f)=(x_{c}x_{e}x_{a},x_{f}x_{b}x_{d},x_{c}x_{f}^{-1},x_{d}^{-1}x_{a},x_{e},x_{f})\]
This change of variables preserves the Haar measure on $\U(N)^{6}$.\footnote{This is because  the normalised Haar measure on $\U(N)^{n}$, or on $G^{n}$ for any compact topological group $G$, is pushed forward onto itself by each of the elementary maps
\begin{itemize}[itemsep=0pt]
\item $(x_{1},x_{2},\ldots,x_{n})\mapsto (x_{1}^{-1},x_{2},\ldots,x_{n})$
\item $(x_{1},x_{2},\ldots,x_{n})\mapsto (x_{1}x_{2},x_{2},\ldots,x_{n})$
\item $(x_{1},\ldots,x_{n})\mapsto (x_{\sigma(1)},\ldots,x_{\sigma(n)})$, where $\sigma$ is any permutation of $\{1,\ldots,n\}$
\end{itemize}
and it is not difficult to check that our change of variables can be obtained as a composition of these maps. 

Interestingly, these elementary operations are exactly the Nielsen transformations, which generate the group of automorphisms of the free group of rank $n$ (see \cite{LyndonSchupp}). Thus, the random homomorphism from the free group $\mathbb F_{n}$ to a compact topological group $G$ constructed by picking a basis of $\mathbb F_{n}$ and sending this basis to a uniformly chosen element of~$G^{n}$ does not depend, in distribution, on the basis of $\mathbb F_{n}$ used to construct it. In particular, the distribution of the image of every element of the free group is intrinsically defined, and one may for instance wonder, for specific or for general $G$, which elements of $\mathbb F_{n}$ are sent to a uniformly distributed element of $G$. I am grateful to the referee for pointing out to me that this problem was solved for finite groups in \cite{PuderParzanchevski}.
} Thus, we find
\[\E[f(H_{\ell})]=\int_{\U(N)^{4}} f(g^{-1}yh^{-1}gzh) p_{s}(y)p_{t}(z)p_{u}(g)p_{v}(h)\d g \d h \d y \d z\]
after integrating with respect to $e$ and $f$ which do not appear in the integrand. Thus, considering four independent Brownian motions $G,H,Z,Y$ on $\U(N)$, we find the equality in distribution
\begin{equation}\label{eq:H}
H_{\ell} \build{=}_{}^{\text{dist.}} G_{u}^{-1}Y_{s}H_{v}^{-1}G_{u}Z_{t}H_{v}
\end{equation}

The quantity $\E[\tr(H_{\ell})]$ appears now as a function of the four real parameters $s,t,u,v$ and we can use stochastic calculus to differentiate it with respect to each of them. In fact, using the first assertion of \eqref{eq:tr}, which in the language of Brownian motion reads $\E[Y_{s}]=e^{-\frac{s}{2}}I_{N}$ and $\E[Z_{t}]=e^{-\frac{t}{2}}I_{N}$, we can simplify the problem to
\[\E[\tr(H_{\ell})]=e^{-\frac{s+t}{2}}\E[\tr(G_{u}^{-1}H_{v}^{-1}G_{u}H_{v})]\]
The expectation in the right-hand side of this equality is a symmetric function of $u$ and $v$. Using stochastic calculus, we find
\begin{equation}\label{eq:DEyang}
\partial_{u}\E[\tr(G_{u}^{-1}H_{v}^{-1}G_{u}H_{v})]=-\E[\tr(G_{u}^{-1}H_{v}^{-1}G_{u}H_{v})]+\E[\tr(H_{v}^{-1})\tr(H_{v})]
\end{equation}
The new function $\E[\tr(H_{v}^{-1})\tr(H_{v})]$ of $v$ can in turn be computed using It\={o}'s formula, since it is equal to $1$ when $v=0$ and satisfies the differential equation
\[\partial_{v}\E[\tr(H_{v}^{-1})\tr(H_{v})]=-\E[\tr(H_{v}^{-1})\tr(H_{v})]+\frac{1}{N^{2}}\]
which is solved in 
\begin{equation}\label{eq:trbrbr}
\E[\tr(H_{v}^{-1})\tr(H_{v})]=\frac{1}{N^{2}}(1-e^{-v})+e^{-v}
\end{equation}
Replacing in \eqref{eq:DEyang} and solving, we find finally
\begin{empheq}[box=\mybox]{equation} \label{eq:tryang}
\E[\tr(H_{\ell})]=e^{-\frac{s+t}{2}} \Big(e^{-u}+e^{-v}-e^{-(u+v)}+\frac{1}{N^{2}}(1-e^{-u})(1-e^{-v})\Big)
\end{empheq}
and, letting $N$ tend to infinity,
\begin{empheq}[box=\mybox]{equation} \label{eq:tryangMF}
\lim_{N\to\infty}\E[\tr(H_{\ell})]=e^{-\frac{s+t}{2}} \big(e^{-u}+e^{-v}-e^{-(u+v)}\big)
\end{empheq}

\begin{figure}[h!]
\begin{center}
\includegraphics{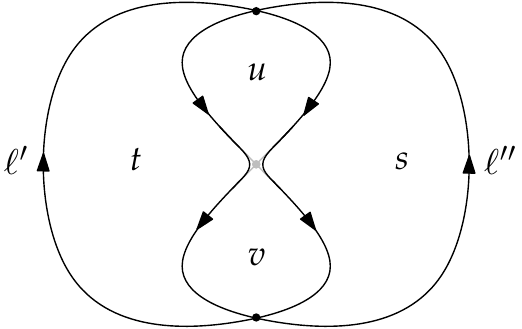}
\caption{\label{fig:desingeight} \small The loops $\ell'$ and $\ell''$ are obtained from $\ell$ by an operation that will feature prominently in Section \ref{sec:MMeqns}.
}
\end{center}
\end{figure}

We did these computations without taking great care of a possible geometric meaning of the successive steps. Anticipating our discussion of the Makeenko--Migdal equations, it is interesting to check that 
\begin{equation}\label{eq:MM8}
(\partial_{u}+\partial_{v}-\partial_{s}-\partial_{t})\E[\tr(H_{\ell})]=e^{-\frac{s+t}{2}}\Big(e^{-(u+v)}+\frac{1}{N^{2}}(1-e^{-(u+v)})\Big)=\E[\tr(H_{\ell'})\tr(H_{\ell''})]
\end{equation}
where $\ell'$ and $\ell''$ are the loops drawn on Fig. \ref{fig:desingeight}.

Perhaps even more interesting than the fact that \eqref{eq:MM8} holds, which after all is a consequence of Theorem \ref{thm:MME}, is the observation that \eqref{eq:MM8} does not seem to be easily guessed from \eqref{eq:H} and It\={o}'s formula. More precisely, It\={o}'s formula allows us to give an expression of the left-hand side of \eqref{eq:MM8} and it is not obvious that this expression coincides with the right-hand side of \eqref{eq:MM8}. We take this as a sign that the Makeenko--Migdal equations give an information that is practically non trivial.

\subsection{The case of the sphere: a not so simple loop}

Computations involving the Yang--Mills holonomy process on the sphere, although in principle based on the same formulas as in the case of the plane, are in general much more complicated. This can be explained by the fact that, as we indicated in Section \ref{sec:structYMHP}, the stochastic core of the Yang--Mills holonomy process on a sphere is a Brownian {\em bridge} on $\U(N)$, or on the compact Lie group $G$, instead of a Brownian motion.

In this section, we are going to illustrate some of the difficulties that one meets when working on a sphere. The first is that the partition function is not equal to $1$ anymore. Instead, according to \eqref{prop:Z}, it is given, on a sphere of total area $T$, by
\[Z_{S^{2}}(T)=p_{T}(I_{N})=\|p_{\frac{T}{2}}\|^{2}_{L^{2}(\U(N))}=\sum_{\alpha\in \widehat \U(N)} e^{-\frac{T}{2} c_{2}(\alpha)} d_{\alpha}^{2}\]
This is also an expression in which nothing is specific to $\U(N)$: it is valid for any compact Lie group.\footnote{Note that $T$, which used to denote the coupling constant in \eqref{prop:Z}, now denotes the total area of our surface. This is not a problem because the only meaningful quantity is the product of the coupling constant by the total area of the surface.}

The most basic question about the Yang--Mills holonomy process on the sphere is the analogue to the question that we treated in Section \ref{sec:simpleplane}, namely to compute the expectation of the normalised trace of the holonomy along a simple loop $\ell$ enclosing a domain of area $t$. The Driver--Sengupta formula yields the following expression for this expectation:
\begin{equation}\label{eq:expsph}
\E[\tr(H_{\ell})]=\frac{1}{Z_{S^{2}}(T)} \int_{\U(N)} \tr(x) p_{t}(x)p_{T-t}(x^{-1})\d x
\end{equation}
Using the Fourier expansion of the heat kernel, one finds 
\[\E[\tr(H_{\ell})]=\frac{1}{Z_{S^{2}}(T)}
\sum_{\lambda,\mu\in \widehat\U(N)} e^{-c_{2}(\lambda)\frac{t}{2}-c_{2}(\mu)\frac{T-t}{2}} d_{\lambda} d_{\mu} \int_{\U(N)} \tr(x)\chi_{\lambda}(x)\chi_{\mu}(x^{-1})\d x\]
The integral can be computed thanks to Pieri's rule: it is equal to $0$ unless $\mu$ is obtained from $\lambda$ by adding $1$ to exactly one component, in which case it is equal to $1$. We write $\lambda \nearrow \mu$ when this happens. Thus,
\begin{equation}\label{eq:harmsph}
\E[\tr(H_{\ell})]=\frac{1}{Z_{S^{2}}(T)}
\sum_{\lambda\in \widehat\U(N)} e^{-c_{2}(\lambda)\frac{T}{2}} d_{\lambda}^{2} \underbrace{\bigg[\sum_{\substack{\mu\in \widehat \U(N)\\ \lambda \nearrow \mu}}
e^{-(c_{2}(\mu)-c_{2}(\lambda))\frac{T-t}{2}} \frac{d_{\mu}}{d_{\lambda}}\bigg]}_{f_{1}(\lambda)} 
\end{equation}
It seems difficult to give an expression of $\E[\tr(H_{\ell})]$ much simpler than \eqref{eq:expsph} or \eqref{eq:harmsph} which, as is hardly necessary to emphasize, is much more complicated than the one that we obtained in the case of the plane.\footnote{Let us drive the point home: \eqref{eq:harmsph}, once made fully explicit using \eqref{eq:dimcas}, is the exact analogue of the $e^{-\frac{t}{2}}$ that we see in the second assertion of \eqref{eq:tr}.}

It is, however, possible to analyse the limit of this quantity as $N$ tends to infinity. A first step in this direction is based on the realisation that Pieri's rule is simple, and the quantity between square brackets, which we denote by $f_{1}(\lambda)$ is a finite sum and can be written explicitly using \eqref{eq:dimcas}:
\[f_{1}(\lambda)=e^{-\frac{T-t}{2}}\sum_{i=1}^{N} {\mathbbm 1}_{\{i=1 \text{ or } \lambda_{i-1}>\lambda_{i}\}} e^{-(T-t)\big(\lambda_{i}+\frac{N-2i+1}{2}\big)}\prod_{\substack{1\leq j \leq N \\ j\neq i}} \Big(1+\frac{1}{\lambda_{i}-\lambda_{j}+j-i}\Big)\]
This suggest to associate to the highest weight $\lambda$ the decreasing sequence $l=(l_{1}>\ldots >l_{N})$ of half-integers defined by
\[l_{i}=\lambda_{i}+\frac{N-2i+1}{2}\]
so that
\[f_{1}(\lambda)=e^{-\frac{T-t}{2}}\sum_{i=1}^{N} {\mathbbm 1}_{\{i=1 \text{ or } \lambda_{i-1}>\lambda_{i}\}} e^{-\frac{T-t}{N}l_{i}} \prod_{\substack{1\leq j \leq N \\ j\neq i}} \Big(1+\frac{1}{l_{i}-l_{j}}\Big)\]

Let us now introduce the probability measure $\pi_{N,T}$ on $\widehat \U(N)$ such that for every highest weight $\lambda$, one has
\[\pi_{N,T}(\{\lambda\})\propto e^{-c_{2}(\lambda)\frac{T}{2}}d_{\lambda}^{2}\]
 Then \eqref{eq:harmsph} can be written more compactly as
\begin{equation}\label{eq:intUhat}
\E[\tr(H_{\ell})]=\int_{\widehat \U(N)} f_{1}(\lambda) \d \pi_{N,T}(\lambda)
\end{equation}
Moreover, there exists for each integer $n\geq 2$ a function $f_{n}$ on $\widehat\U(N)$, not very different from $f_{1}$, and the integral of which against $\pi_{N,T}$ yields $\E[\tr(H_{\ell}^{n})]$.

We would like to express that, as $N$ tends to infinity, the measure $\pi_{N,T}$ concentrates on a few highest weights, characterised by a certain limiting shape. One unpleasant feature of \eqref{eq:intUhat} in this respect is that the set on which the integral is taken, namely $\widehat\U(N)$, depends on $N$. It is thus uneasy to formulate a concentration result. One classical and efficient way around this problem is to associate to each highest weight $\lambda$ its {\em empirical measure}
\[\hat\mu_{\lambda}=\frac{1}{N}\sum_{i=1}^{N} \delta_{\frac{l_{i}}{N}}=\frac{1}{N}\sum_{i=1}^{N} \delta_{\frac{1}{N}(\lambda_{i}+\frac{N-2i+1}{2})}\]

\begin{figure}[h!]
\begin{center}
\includegraphics{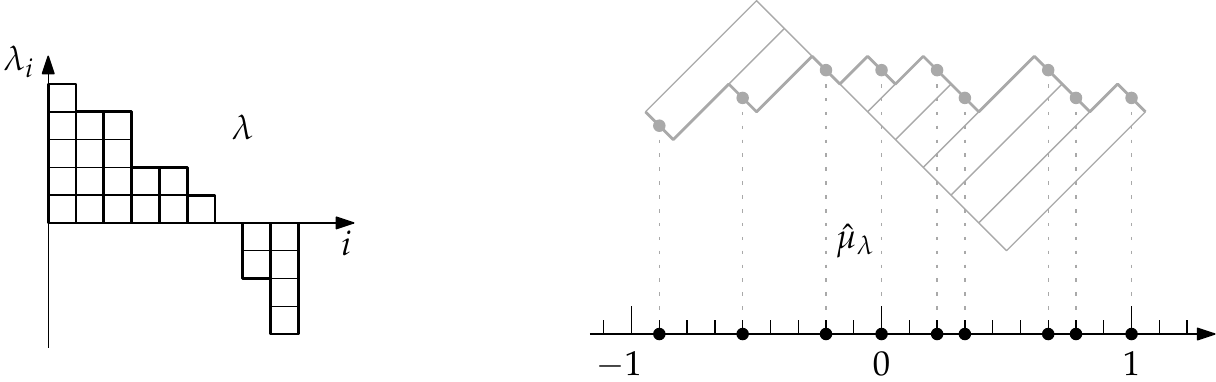}
\caption{\label{fig:lambdaell} \small With $N=9$, the highest weight $\lambda=(5,4,4,2,2,1,0,-2,-4)$ drawn in the style of a Young diagram, and its empirical measure. Each dot represents $\frac{1}{9}$ of mass and any two dots are distant by a multiple of $\frac{1}{9}$.}
\end{center}
\end{figure}

Pushing the probability measure $\pi_{N,T}$ forward by the map $\lambda \mapsto \hat\mu_{\lambda}$ yields a probability measure, which we denote by $\Pi_{N,T}$, on the set of probability measures on the real line. It is possible to predict the behaviour of this probability as $N$ tends to infinity by writing $c_{2}(\lambda)$ and $\d_{\lambda}$ in terms of the empirical measure of $\lambda$. Up to some inessential terms (see \cite[Eq. (24)]{LevyMaida} for complete expressions), one finds
\[c_{2}(\lambda)\simeq N^{2} \int_{\R} x^{2} \d\hat\mu_{\lambda}(x) \ \ \text{ and } \ \ d_{\lambda}^{2} \simeq \exp \bigg[- N^{2} \int_{\{(x,y)\in \R^{2}, x\neq y\}} -\log |x-y| \d\hat\mu_{\lambda}(x) \d\hat\mu_{\lambda}(y)\bigg]\]
Introducing, for every probability measure $\mu$, the quantity
\[\mathcal J_{T}(\mu)=\int_{\{(x,y)\in \R^{2}, x\neq y\}} -\log |x-y| \d\mu(x) \d\mu(y)+\frac{T}{2}\int_{\R} x^{2}\d\mu(x)\]
we see that the probability measure $\Pi_{N,T}$ assigns to any probability measure $\mu$ that is the empirical measure of a highest weight a mass proportional to
\[\Pi_{N,T}(\{\mu\})\propto  \exp(-N^{2}\mathcal J_{T}(\mu))\]
In the large $N$ limit, it seems plausible that $\Pi_{N,T}$ will concentrate on the minimisers, or even better, on the unique minimiser of the functional $\mathcal J_{T}$. This turns out to be true, with a little twist that we will explain and contributes to making the story much more interesting than it already is. Let us summarise the main results on which one can ground a rigorous analysis of the situation.

\begin{figure}[h!]
\begin{center}
\includegraphics{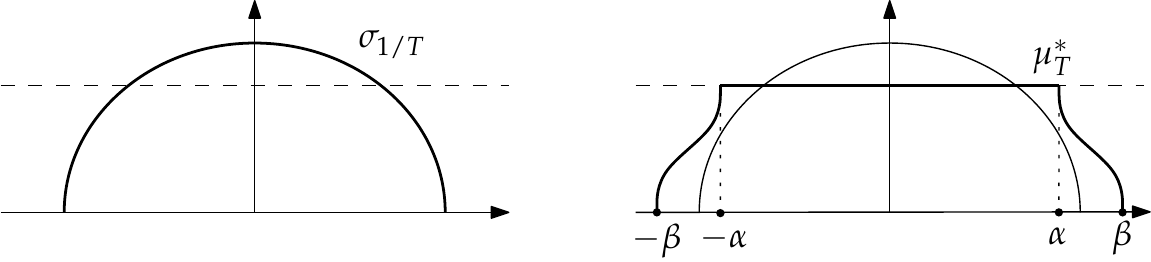}
\caption{\label{fig:sigmamu} \small For $T>\pi^{2}$, the absolute minimiser of the functional $\mathcal J_{T}$ does not belong to the class of probabilities on $\R$ with a density not exceeding $1$. The minimiser within this class is represented on the right. Its density is identically equal to $1$ on an interval in the middle of its support, and given by elliptic functions outside this interval.}
\end{center}
\end{figure}

\begin{itemize}
\item Minimising the functional $\mathcal J_{T}$ on the space of all probability measures on $\R$ is one of the simplest examples of a rich and well-developed theory which is, for example, exposed in the book of Edward Saff and Vilmos Totik \cite{SaffTotik}. This is also a very common problem in random matrix theory. Indeed, the unique minimiser of $\mathcal J_{T}$ is Wigner's semi-circular distribution with variance $\frac{1}{T}$:
\begin{equation}\label{eq:minsc}
\d\sigma_{1/T}(x)=\frac{T}{2\pi}\sqrt{\frac{4}{T}-x^{2}}\  \mathbbm 1_{\big[-\frac{2}{\sqrt{T}},\frac{2}{\sqrt{T}}\big]}(t)\d t
\end{equation}
\item The fact that the measure $\Pi_{N,T}$ concentrates, as $N$ tends to infinity, to the minimiser of $\mathcal J_{T}$ is a special case of a principle of large deviations proved by Alice Guionnet and Myl\`ene Ma\"\i da in \cite{GuionnetMaida}. However, the minimiser of $\mathcal J_{T}$ that one must consider is not the absolute minimiser on the set of all probability measures on $\R$. Indeed, for all $N\geq 1$, the measure $\Pi_{N,T}$ is supported by the set of empirical measures of highest weights of $\U(N)$, which form a rather special set of probability measures. A distinctive feature of these measures is that they are atomic, with atoms of mass $\frac{1}{N}$ spaced by integer multiples of $\frac{1}{N}$. Weak limits, as $N$ tends to infinity, of such measures can only be absolutely continuous with respect to the Lebesgue measure on $\R$, with a density not exceeding~$1$: a class of probability measures that we will denote by $\mathcal L(\R)$. The result of Guionnet and Ma\"\i da asserts that the measure $\Pi_{N,T}$ concentrates exponentially fast, as $N$ tends to infinity, around the unique minimiser $\mu_{T}^{*}$ of $\mathcal J_{T}$ on the closed set $\mathcal L(\R)$.

\item The problem of minimising $\mathcal J_{T}$ under the constraint of having a density not exceeding~$1$ is a problem that is, in principle, just as well understood as the unconstrained problem. The book \cite{SaffTotik} contains results ensuring the existence and uniqueness of the minimiser, and others allowing one to determine its support. In fact, the measure $\sigma_{1/T}$ given by \eqref{eq:minsc}, and which is the absolute minimiser of $\mathcal J_{T}$, is absolutely continuous with a maximal density of $\sqrt{T}/\pi$, so that it belongs to $\mathcal L(\R)$ provided $T\leq \pi^{2}$. For $T>\pi^{2}$, the constraint becomes truly restrictive, and one must make do with a probability measure which is, in $\mathcal L(\R)$, the best available substitute for $\sigma_{1/T}$. The actual determination of this minimiser $\mu^{*}_{T}$ is, depending on one's background, a more or less elementary exercise in Riemann--Hilbert theory, and involves manipulating elliptic functions. The density of~$\mu_{T}^{*}$ for $T>\pi^{2}$ is represented on Fig. \ref{fig:sigmamu}. An exact expression of this density can be found in \cite[Eq. (37)]{LevyMaida}.
\end{itemize}

Having established the exponential concentration, as $N$ tends to infinity, of the measure $\Pi_{N,T}$ around $\mu^{*}_{T}$, it is possible to come back to our initial problem of computing $\E[\tr(H_{\ell})]$. After noticing that $f_{1}(\lambda)$ can be expressed as a functional $F_{1}(\hat\mu_{\lambda})$ of the empirical measure of $\lambda$, it can be guessed that $\E[\tr(H_{\ell})]$ is related to $F_{1}(\mu^{*}_{T})$. Antoine Dahlqvist and James Norris were the first to rigorously and successfully pursue this line of reasoning, and to obtain the following remarkably elegant result.

\begin{theorem}[Dahlqvist--Norris \cite{DahlqvistNorris}] Let $\rho_{T}$ denote the density of the minimiser $\mu^{*}_{T}$. Then, for all integer $n\geq 0$, one has
\begin{equation}\label{eq:momsphere}
\lim_{N\to \infty}\E[\tr(H_{\ell}^{n})]=\lim_{N\to \infty}\E[\tr(H_{\ell}^{-n})]=\frac{1}{n\pi} \int_{\R} \cosh \Big(\frac{nx}{2}(T-2t)\Big) \sin (n\pi \rho_{T}(x))\d x
\end{equation}
\end{theorem}

To conclude this long discussion of the simple loop on the sphere, let us mention another result for the statement of which we have all the concepts at hand. Our description of the behaviour of the measure $\Pi_{N,T}$ suggests that the partition function itself is dominated by the contribution of the highest weights that have an empirical measure close to $\mu^{*}_{T}$. This is indeed true, and the fact that the shape of $\mu^{*}_{T}$ changes suddenly when $T$ crosses the critical value $\pi^{2}$ gives rise to a phase transition, in this case of third order, first discovered by Douglas and Kazakov, and named after them. It was first proved rigourously, in a slightly different but equivalent language, by Karl Liechty and Dong Wang in \cite{LiechtyWang}, and by Myl\`ene Ma\"\i da and the author in \cite{LevyMaida}.

\begin{theorem}[Douglas--Kazakov phase transition] The free energy of the Yang--Mills model on a sphere of total area $T$ is given by
\[F(T)=\lim_{N\to \infty} \frac{1}{N^{2}}\log Z_{S^{2}}(T)=\frac{T}{24}+\frac{3}{2}-\mathcal J_{T}(\mu^{*}_{T})\]
The function $F$ is of class $C^{2}$ on $(0,\infty)$ and smooth on $(0,\infty)\setminus\{\pi^{2}\}$. The third derivative of $F$ admits a jump discontinuity at $\pi^{2}$.
\end{theorem}

\begin{figure}[h!]
\begin{center}
\includegraphics[width=7cm]{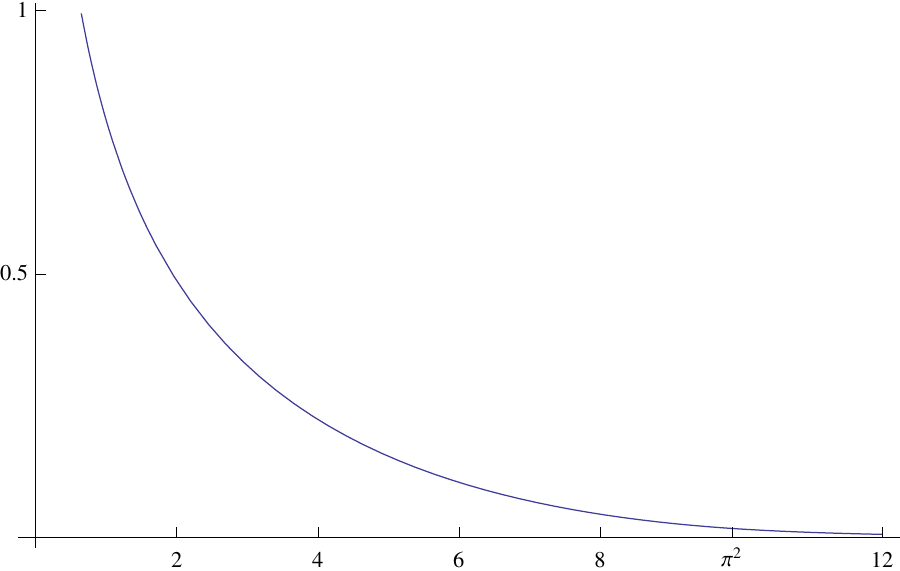}\hspace{1cm}\includegraphics[width=7cm]{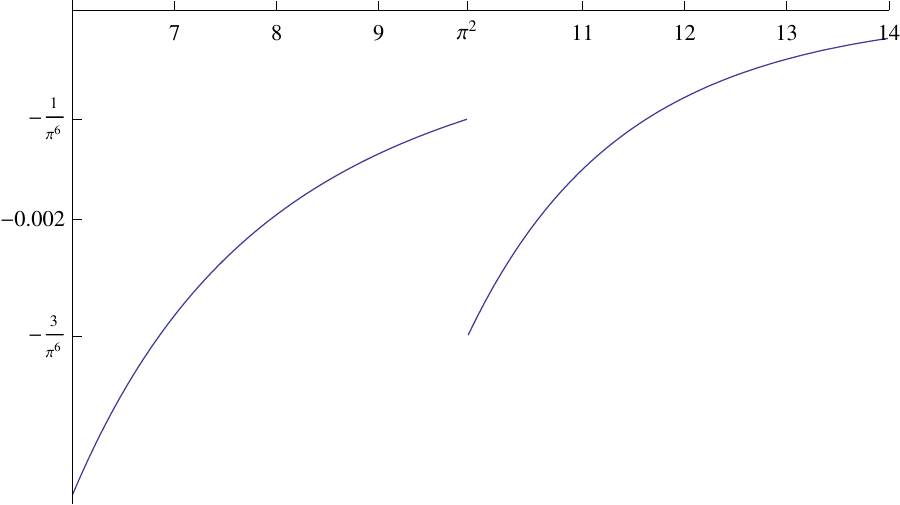}
\caption{\label{fig:GF} \small The graphs of $T\mapsto F(T)$ (on the left) and of $T\mapsto F^{(3)}(T)$ near $T=\pi^{2}$ (on the right).}
\end{center}
\end{figure}

This phase transition is not one that is easily detected numerically, as Fig. \ref{fig:GF} shows.

\section{The Makeenko--Migdal equations}\label{sec:MMeqns}
\subsection{First approach}\label{sec:faMM}

\begin{figure}[h!]
\begin{center}
\includegraphics[width=12cm]{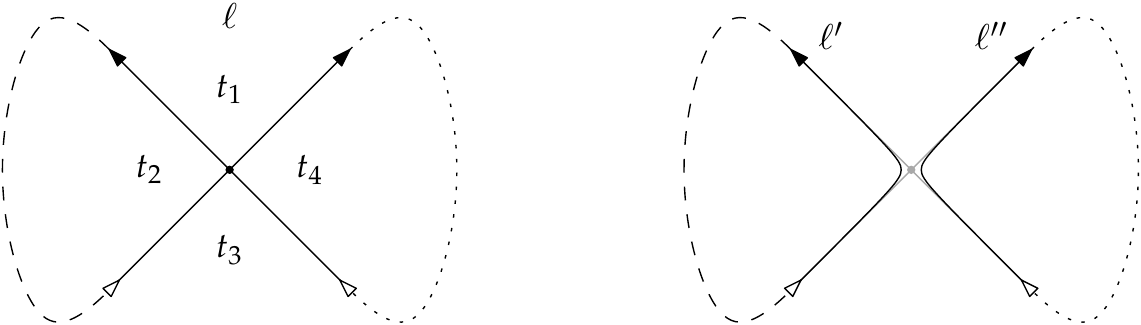}
\caption{\label{fig:desing} \small On the left, we see a loop $\ell$ around a generic self-intersection point. The dotted and dashed part of $\ell$ can be arbitrarily complicated, and can meet many times outside the small region of the surface that we are focusing on. It is nevertheless true that after escaping  this small region through the North-East corner (resp. North-West corner), the first time $\ell$ comes back is through the South-East corner (resp. South-West corner). This is why the `desingularisation' operation illustrated on the right produces exactly two loops, that we call $\ell'$ and $\ell''$.   }
\end{center}
\end{figure}

It is now time that we discuss the equations discovered by Yuri Makeenko and Alexander Migdal and which give their title to these notes. These equations are a powerful tool for the study of the Wilson loop expectations of which we gave a few examples in the previous section. They are related to the approach that we called dynamical, in which an expectation of the form $\E[\tr(H_{\ell})]$, where $\ell$ is some nice loop on a surface $M$, is seen as a function of the areas of the faces cut by $\ell$ on the surface $M$. The Makeenko--Migdal equations give a remarkably  elegant expression of the alternated sum of the derivatives of $\E[\tr(H_{\ell})]$ with respect to the areas of the four faces that surround a generic point of self-intersection of $\ell$. This expression is of the form $\E[\tr(H_{\ell'})\tr(H_{\ell''})]$, where $\ell'$ and $\ell''$ are two loops obtained from $\ell$ by a very simple operation at this point of self-intersection $\ell$. This operation consists in taking the two incoming strands of $\ell$ at this point and connecting them with the two outgoing strands in the `other' way, the way that is not realised by $\ell$, see Fig. \ref{fig:desing}.

On this figure, we see four faces around the self-intersection point, which need not be pairwise distinct. We denote their areas by $t_{1},t_{2},t_{3},t_{4}$ as indicated on Fig. \ref{fig:desing}. The Makeenko--Migdal equation in this case reads 
\begin{empheq}[box=\mybox]{equation} \label{eq:MM0}\tag{MM}
\bigg(\frac{\partial}{\partial t_{1}}-\frac{\partial}{\partial t_{2}}+\frac{\partial}{\partial t_{3}}-\frac{\partial}{\partial t_{4}}\bigg)\E[\tr(H_{\ell})]=\E[\tr(H_{\ell'})\tr(H_{\ell''})]
\end{empheq}
The relation \eqref{eq:eyeMM1}, that we derived earlier in an elementary way, is an instance of this equation. 

The relation \eqref{eq:MM0} would become particularly useful if we could combine it with a result saying that $\E[\tr(H_{\ell'})\tr(H_{\ell''})]=\E[\tr(H_{\ell'})]\E[\tr(H_{\ell''})]$. A crucial fact is that this equality, which is of course false  in general, becomes true in the large $N$ limit in all cases where this limit has been studied, that is, on the plane and on the sphere. It corresponds to a concentration phenomenon, namely to the fact that the complex-valued random variable $\tr(H_{\ell})$ converges, in the large $N$ limit, to a deterministic complex, indeed real number $\Phi(\ell)$. This behaviour is expected to occur on any compact surface, and the function $\Phi:\Loop(M)\to \R$, whose existence has so far been proved when $M$ is the plane or the sphere, is called the {\em master field}.

In the large $N$ limit, the Makeenko--Migdal equation \eqref{eq:MM0} becomes a kind of differential equation satisfied by this master field $\Phi$:
\begin{empheq}[box=\mybox]{equation} \label{eq:MMinf}\tag{$\text{MM}_{\infty}$}
\bigg(\frac{\partial}{\partial t_{1}}-\frac{\partial}{\partial t_{2}}+\frac{\partial}{\partial t_{3}}-\frac{\partial}{\partial t_{4}}\bigg)\Phi(\ell)=\Phi(\ell')\Phi(\ell'')
\end{empheq}

On the plane, we will see that this equation, together with the very simple equation \eqref{eq:tr}, essentially characterises the function $\Phi$.

\subsection{Makeenko and Migdal's proof}\label{sec:MMproof}

Makeenko and Migdal discovered the relation \eqref{eq:MM0}, and the extensions that we will describe later, by doing a very clever integration by parts in the functional integral with respect to the Yang--Mills measure (see \eqref{eq:YM}) that defines a Wilson loop expectation:
\[\E[\tr(H_{\ell})]=\frac{1}{Z}\int_{\A} \tr(\hol(\omega,\ell)) e^{-\frac{1}{2}S_{\YM}(\omega)}\d \omega
\]
or instead, as we will explain, in a closely related integral (see \cite{MakeenkoMigdal} and \cite{DubinMakeenko}). That this integration by parts  performed in an ill-defined integral yields as a final product a perfectly meaningful formula, makes Makeenko and Migdal's original derivation the more intriguing. It is described in mathematical language in the introduction of \cite{LevyMF}, but this derivation is so beautiful that we reproduce its description here.

The finite-dimensional prototype of the so-called Schwinger--Dyson equations, obtained by integration by parts in functional integrals, is the fact that for all smooth function $f:\R^{n}\to \R$ with bounded differential, and for all $h\in \R^{n}$, the equality
\[\int_{\R^{n}} d_{x}f(h) e^{-\frac{1}{2}\|x\|^{2}}\; dx=\int_{\R^{n}}\langle x,h\rangle f(x)e^{-\frac{1}{2}\|x\|^{2}}\; dx\]
holds. This equality ultimately relies on the invariance by translation of the Lebesgue measure on $\R^{n}$ and it can be proved by writing
\[0=\frac{d}{dt}_{|t=0}\int_{\R^{n}}f(x+th)e^{-\frac{1}{2}\|x+th\|^{2}}\; dx\]

In our description of the Yang--Mills measure $\mu_{\YM}$ (see \eqref{eq:YM}), we mentioned that the measure $\d \omega$ on the space $\A$ of connections was meant to be a kind of Lebesgue measure, invariant by translations. This is the key to the derivation of the Schwinger--Dyson equations, as we will now explain. In what follows, we will use the differential geometric language introduced in Section \ref{sec:YMconn}.

Let $\psi:\A\to \R$ be an observable, that is, a function. In general, we are interested in the integral of $\psi$ with respect to the measure $\mu_{\YM}$. The tangent space to the affine space $\A$ is the linear space $\Omega^{1}(M)\otimes \Ad(P)$. To say that the measure $\d \omega$ is translation invariant means that for every element $\eta$ of this linear space,
\[0=\frac{d}{dt}_{|t=0}\int_{\A}\psi(\omega+t\eta)e^{-\frac{1}{2}S_{\YM}(\omega+t\eta)}\d\omega\]
and the Schwinger--Dyson equations follow in their abstract form
\begin{equation}\label{SD intro}
\int_{\A}d_{\omega}\psi(\eta) \d\mu_{\YM}(\omega)=\frac{1}{2}\int_{\A}\psi(\omega)d_{\omega}S_{\YM}(\eta)\d\mu_{\YM}(\omega)
\end{equation}
The directional differential of the Yang--Mills action is well known (see for example \cite{Bleecker}) and most easily expressed using the covariant exterior differential $d^{\omega}:\Omega^{0}(M)\otimes \Ad(P)\to \Omega^{1}(M)\otimes \Ad(P)$ defined by $d^{\omega}\alpha=d\alpha+[\omega\wedge \alpha]$. It is given by
\[d_{\omega}S_{\YM}(\eta)=2\int_{M}\langle \eta \wedge d^{\omega} *\!\Omega\rangle\]
The problem is now to apply this formula to a well-chosen observable $\psi$ and to differentiate in the right direction.

Given a loop $\ell$ on $M$, Makeenko and Migdal applied \eqref{SD intro} to the observable defined by choosing a skew-Hermitian matrix $X\in \u(N)$ and setting, for all $\omega\in \A$,
\begin{equation}\label{eq:psiX1}
\psi_{X}(\omega)=\Tr(X\,\hol(\omega,\ell))
\end{equation}
To make this definition perfectly meaningful, one needs to choose a reference point in the fibre of $P$ over the base point of $\ell$: we will assume that such a point has been chosen and fixed, and compute holonomies with respect to this point. 

Let us choose a parametrisation $\ell:[0,1]\to M$ of $\ell$. The directional derivative of the observable $\psi_{X}$ in the direction of a $1$-form $\eta\in \Omega^{1}(M)\otimes \Ad(P)$ is given by
\begin{equation}\label{dpsiomega}
d_{\omega}\psi_{X}(\eta)=-\int_{0}^{1}\Tr\left(X\,\hol(\omega,\ell_{[s,1]})\eta(\dot \ell(s)) \hol(\omega,\ell_{[0,s]})\right) \; ds
\end{equation}
where we denote by $\ell_{[a,b]}$ the restriction of $\ell$ to the interval $[a,b]$.\footnote{At first glance, \eqref{dpsiomega} may seem to require the choice of a point in $P_{\ell(s)}$ for each $s$, but in fact it does not, for the way in which the two holonomies and the term $\eta(\dot \ell(s))$ would depend on the choice of this point cancel exactly.}

One must now choose the direction of differentiation $\eta$. Let us assume that $\ell$ is a nice loop which around each point of self-intersection looks like the left half of Fig. \ref{fig:desing}. Let us assume that for some $s_{0}\in (0,1)$, we have $\ell(s_{0})=\ell(0)$ and $\det(\dot \ell(0),\dot \ell(s_{0}))=1$. Makeenko and Migdal choose for $\eta$ a distributional $1$-form supported at the self-intersection point $\ell(0)$, which one could write as\footnote{It may seem that we are progressively letting go of the intrinsic character of our construction, but the interested reader can check that everything is still geometrically meaningful at this point.}
\[\forall m\in M, \forall v\in T_{m}M, \; \eta_{m}(v)=\delta_{m,\ell(0)} \det(\dot \ell(0),v) X\]
with $\det(\dot\ell(0),v)$ denoting the determinant of the two vectors $\dot\ell(0)$ and $v$. With this choice of $\eta$, the directional derivative of $\psi_{X}$ is given by
\begin{equation}\label{eq:psiX2}
d_{\omega}\psi_{X}(\eta)=-\Tr\left(X\,\hol(\omega,\ell_{[s_{0},1]}) X\, \hol(\omega,\ell_{[0,s_{0}]})\right)=-\Tr\left(X\,\hol(\omega,\ell') X\, \hol(\omega,\ell'')\right)
\end{equation}
where $\ell'$ and $\ell''$ are the loops defined on the right of Fig. \ref{fig:desing}. 
Recall that $\u(N)$ is endowed with the invariant scalar product $\langle X,Y\rangle=-N\Tr(XY)$. The directional derivative of the Yang--Mills action is thus given by 
\[d_{\omega}S_{\YM}(\eta)=-2\langle X, (d^{\omega}\!*\!\Omega)(\dot \ell(0))\rangle=-2N\Tr\left(X d^{\omega}\!*\!\Omega(\dot \ell(0))\right)\]
or so it seems from a naive computation. We shall soon see that this expression needs to be reconsidered. For the time being, our Schwinger--Dyson equation reads
\begin{equation}\label{eq:SDX}\tag{$\text{SD}_{X}$}
\int_{\A}\Tr\left(X\,\hol(\omega,\ell') X\, \hol(\omega,\ell'')\right)\d\mu_{\YM}(\omega)= 
 N\int_{\A}\Tr(X\,\hol(\omega,\ell)) \Tr(X\, d^{\omega}\! * \! \Omega(\dot \ell(0))) \d\mu_{\YM}(\omega)
\end{equation}
Let us add the equalities \eqref{eq:SDX} obtained by letting $X$ take all the values $X_{1},\ldots,X_{N^{2}}$ of an orthonormal basis of $\u(N)$. With the scalar product which we chose, the relations\footnote{These relations are strictly equivalent to \eqref{eq:KK}. They are, in one form or the other, {\em the} fundamental fact of all this story.}
\begin{empheq}[box=\mybox]{equation} \label{eq:magic}
\sum_{k=1}^{N^{2}}\Tr(X_{k}AX_{k}B)=-\frac{1}{N}\Tr(A)\Tr(B) \mbox{ and } \sum_{k=1}^{N^{2}}\Tr(X_{k}A)\Tr(X_{k}B)=-\frac{1}{N}\Tr(AB)
\end{empheq}

hold for any two matrices $A$ and $B$, so that we find
\begin{equation*}
\int_{\A}\tr\left(\hol(\omega,\ell'))\tr(\hol(\omega,\ell'')\right)\d\mu_{\YM}(\omega)=
\int_{\A}\tr\left(\hol(\omega,\ell) d^{\omega}\! * \! \Omega(\dot \ell(0))\right)\d \mu_{\YM}(\omega).
\end{equation*}

The left-hand side of this equation is the right-hand side of \eqref{eq:MM0}. The last and most delicate heuristic step is to interpret the right-hand side of this equation. For this, we must understand the term $d^{\omega}*\!\Omega(\dot \ell(0))$ and we do this by combining two facts: the fact that $d^{\omega}$ acts by differentiation in the horizontal direction and the fact that $*\Omega$ computes the holonomy along infinitesimal rectangles. We must also remember that this term comes from the computation of the exterior product of the distributional form $\eta$ with the form $d^{\omega}*\!\Omega$. It turns out that, instead of a derivative in the horizontal direction with respect to $s$ at $s=0$, we should think of the difference between the values at $0^{+}$ and at $0^{-}$, which we denote by $\Delta_{|s=0}$.

With all this preparation and, it must be said, a small leap of faith,
the right-hand side of the Schwinger--Dyson equation can finally be drawn as follows:
\begin{align*}
&\hspace{-2cm}\Delta_{|s=0} \frac{d}{d\epsilon}_{|\epsilon=0} \int_{\A} \tr\,\hol\Bigg(\omega,\; \raisebox{-7mm}{\scalebox{0.5}{\includegraphics{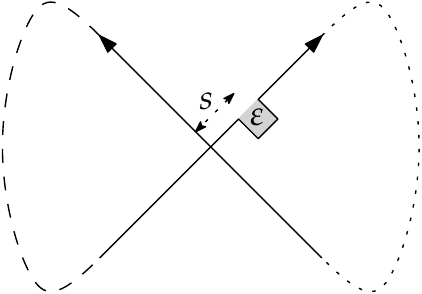}}}\, \Bigg)\d\mu_{YM}(\omega)& \\
\hspace{2cm}= \phantom{-}&\frac{d}{d\epsilon}_{|\epsilon=0} \int_{\A} \tr\,\hol\Bigg(\omega,\;\raisebox{-7mm}{\scalebox{0.5}{\includegraphics{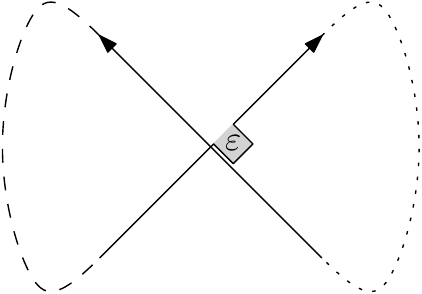}}}\, \Bigg)\d\mu_{YM}(\omega)\\
\hspace{2cm}-&\frac{d}{d\epsilon}_{|\epsilon=0} \int_{\A} \tr\,\hol\Bigg(\omega,\;\raisebox{-7mm}{\scalebox{0.5}{\includegraphics{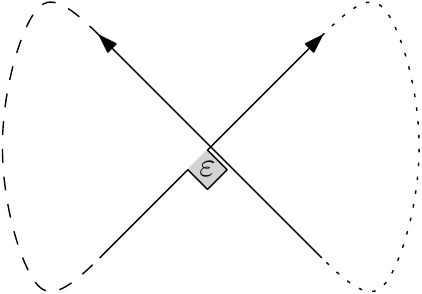}}}\, \Bigg)\d\mu_{YM}(\omega)\\
\end{align*}
This is indeed the left-hand side of the Makeenko--Migdal equation \eqref{eq:MM0}.

\subsection{The equations, their merits and demerits}

\begin{figure}[h!]
\begin{center}
\includegraphics[width=12cm]{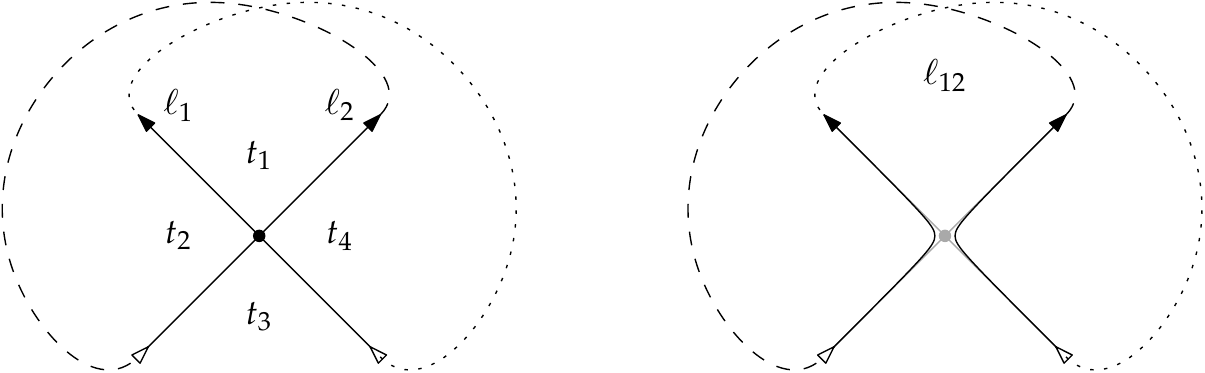}
\caption{\label{fig:desing2} \small When performed at a crossing of two distinct loops $\ell_{1}$ and $\ell_{2}$, the operation of reconnecting the incoming and outgoing strands in the other way that is consistent with orientation produces, from $\ell_{1}$ and $\ell_{2}$, one bigger loop that we denote by $\ell_{12}$. }
\end{center}
\end{figure}

The strategy of proof described in the previous section can be used, and was used by Makeenko and Migdal, to derive equations slightly more general than \eqref{eq:MM0}. Let us indeed consider a collection $\ell_{1},\ldots,\ell_{n}$ of loops on the surface $M$. We assume that these loops are nice and in generic position, in the sense that every crossing between two portions of these loops, be they two portions of the same loop or portions of two different loops, is a simple transverse intersection. Around such a crossing, we see, as before, four faces of the graph cut on $M$ by $\ell_{1},\ldots,\ell_{n}$, and we label the areas of these faces $t_{1},t_{2},t_{3},t_{4}$ as indicated on Fig. \ref{fig:desing} and Fig. \ref{fig:desing2}. The Makeenko--Migdal equations express the alternated sum of the derivatives with respect to $t_{1},t_{2},t_{3},t_{4}$ of $\E[\tr(H_{\ell_{1}})\ldots \tr(H_{\ell_{n}})]$. The equations come in two variants, depending on whether the crossing is between two strands of the same loop (let us call this the case I) or between strands of two distinct loops (the case II). In the case II, where the crossing is between strands of two distinct loops, say $\ell_{1}$ and $\ell_{2}$, the same desingularisation operation explained at the beginning of Section \ref{sec:faMM} gives rise to one new loop $\ell_{12}$, as explained in Fig. \ref{fig:desing2}.

Calling, in all cases, $\ell_{1}$ the loop containing the South-West -- North-East strand, one should replace the observable $\psi_{X}$ defined in \eqref{eq:psiX1} by
\[\psi_{X}(\omega)=\Tr(X\hol(\omega,\ell_{1}))\Tr(\hol(\omega,\ell_{2}))\ldots \Tr(\hol(\omega,\ell_{n}))\]
Then the directional derivative of $\psi_{X}$ is given by
\[d_{\omega}\psi_{X}(\eta)=\left|\!\!\begin{array}{ll}
\Tr\big(X\,\hol(\omega,\ell')X\,\hol(\omega,\ell'')\big)\Tr(\hol(\omega,\ell_{2}))\ldots \Tr(\hol(\omega,\ell_{n})) & \text{(case I)}\\[3pt]
\Tr(X\,\hol(\omega,\ell_{1}))\Tr(X\,\hol(\omega,\ell_{2}))\Tr(\hol(\omega,\ell_{3}))\ldots \Tr(\hol(\omega,\ell_{n}))& \text{(case II)}
\end{array}\right.\]
Then, the key to the computation is, as always, given by the equations \eqref{eq:magic}. The final result, with the current notation, is the following.

\begin{theorem}[Makeenko--Migdal equations]\label{thm:MME} Let $\ell_{1},\ldots,\ell_{n}$ be nice loops on $M$ in generic position. Consider a crossing point of two strands of $\ell_{1}$ (case I) or of one strand of $\ell_{1}$ and one strand of $\ell_{2}$ (case~II). Let $t_{1},t_{2},t_{3},t_{4}$ denote the areas of the four faces around this crossing point, as illustrated on Figs. \ref{fig:desing} and \ref{fig:desing2}. Then, with the notation of these figures,
\begin{empheq}[box=\mybox]{equation*} \label{eq:MMIII}\hspace{-3mm}\bigg(\frac{\partial}{\partial t_{1}}-\frac{\partial}{\partial t_{2}}+\frac{\partial}{\partial t_{3}}-\frac{\partial}{\partial t_{4}}\bigg)\E[\tr(H_{\ell_{1}})\ldots \tr(H_{\ell_{n}})]=\left|\!\!\begin{array}{ll}
\E[\tr(H_{\ell'})\tr(H_{\ell''})\tr(H_{\ell_{2}})\ldots \tr(H_{\ell_{n}})]
 & \text{(I)}\\[5pt]
\frac{1}{N^{2}}\E[\tr(H_{\ell_{12}})\tr(H_{\ell_{3}})\ldots \tr(H_{\ell_{n}})]
& \text{(II)}
\end{array}\right. \hspace{-3mm}
\end{empheq}
It is understood that if two of the four faces around the crossing under consideration are identical, then the corresponding derivative should be taken twice. Moreover, in the case where $M=\R^{2}$, any term corresponding to the derivative with respect to the area of the unbounded face should be ignored.
\end{theorem}

Makeenko and Migdal's original paper on this subject is \cite{MakeenkoMigdal}. The first mathematical proof of the equations was given in \cite{LevyMF}. It was rather long and convoluted, and restricted to the case where the surface $M$ is the plane $\R^{2}$. Three very short and elegant proofs of the equations were then given, still for the case of the plane, by Bruce Driver, Brian Hall and Todd Kemp in \cite{DriverHallKemp}. Immediately after, the same team joined by Franck Gabriel proved  in \cite{DriverGabrielHallKemp} that the equations hold on any compact surface. There is little point in reproducing here the content of these beautiful papers. Let us simply emphasize that the fundamental computations remain those summarised in \eqref{eq:magic}.

In addition to their simplicity, the Makeenko--Migdal equations have one major quality which is the fact that the collection of loops appearing in the right-hand side has one crossing less compared with the original collection of loops. Indeed, the operation of desingularisation replaces the crossing where it takes place by a tangential contact which, to the price of an arbitrarily small deformation of the loops, can be suppressed. This suggests the possibility of a recursive computation of Wilson loop expectations. We will explain in the next section that it is indeed possible to use the Makeenko--Migdal equations to set up a recursive computation of the large $N$ limit of Wilson loop expectations. 

\begin{figure}[h!]
\begin{center}
\includegraphics{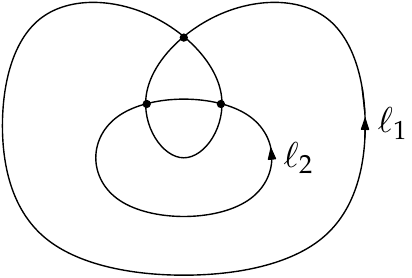}
\caption{\label{fig:CEMM} \small Consider this configuration of two loops on a sphere. It has five faces and three vertices. Moreover, of the three instances of the Makeenko--Migdal equations, two compute the same linear combination of derivatives. There is no hope that the Makeenko--Migdal equations alone will allow one to compute the corresponding Wilson loop expectation.}
\end{center}
\end{figure}

What the Makeenko--Migdal do not do however, is to give a simple formula for the derivative of a Wilson loop expectation with respect to the area of a single face of the graph traced by a given configuration of loops. Only very special linear combinations of these derivatives are accessible. Of course, unless one is working on the plane, the total area of the surface is prescribed and the best one could hope for is a formula describing the variation of the Wilson loop expectations under an arbitrary variation of the areas of the faces that preserves the total area. However, this is, in general, not given by the Makeenko--Migdal equations, see for example Fig. \ref{fig:CEMM}.

It is, in fact, not too difficult to understand what information is available in the Makeenko--Migdal equations. Let us consider $n$ loops $\ell_{1},\ldots,\ell_{n}$ on our surface $M$. Let $F_{1},\ldots,F_{r}$ denote the faces of the graph traced by these loops. Let us identify a vector $(c_{1},\ldots,c_{r})$ of the vector space~$\R^{r}$ with the linear combination of derivatives
\[c_{1}\frac{\partial}{\partial |F_{1}|}+\ldots +c_{r}\frac{\partial}{\partial |F_{r}|}\]
acting on Wilson loop expectations. Let us define the linear subspace $M\subset \R^{r}$ generated by the linear combinations given by the Makeenko--Migdal equations applied at each crossing of the loops $\ell_{1},\ldots,\ell_{n}$. This subspace $M$ is of course contained in the hyperplane $\R^{r}_{0}$ of equation $c_{1}+\ldots +c_{r}=0$. Every element of $\R^{r}$ can naturally be identified with a function on $M$ that is constant on each face of the graph. To each loop $\ell_{i}$, we can associate the unique element ${\sf n}_{\ell_{i}}$ of~$\R^{r}_{0}$ which, as a function on $M$, varies by $1$ across $\ell_{i}$\footnote{A convention must be chosen regarding the definition of a positive crossing of $\ell_{i}$.} and is constant across every other loop. This function is a substitute for the winding number of the loop $\ell_{i}$ on the surface $M$.

It is not difficult to check that it is equivalent, for an element of $\R^{r}$, to be orthogonal, for the simplest scalar product, to the subspace $M$, or to have a constant jump across every loop, the constant possibly depending on the loop. A more formal statement is the following. We denote by $\sf 1$ the vector $(1,\ldots,1)$.

\begin{proposition} In $\R^{r}$, one has the equality of linear subspaces
\[M={\rm Vect}({\sf 1},{\sf n}_{\ell_{1}},\ldots,{\sf n}_{\ell_{n}})^{\perp}\]
In particular, $\dim M=\dim \R^{r}_{0}-n$.
\end{proposition}

The greater the number of loops, the worse the situation. Even with one single loop, we see that all the information about the Wilson loop expectations is not contained in the Makeenko--Migdal equations. 

It is time to turn to a case where things improve drastically, namely the large $N$ limit of the Wilson loop expectations. 

\subsection{The master field on compact surfaces}

We saw in Section \ref{sec:compWLE} that when $G=\U(N)$, Wilson loop expectations tend to take simpler forms in the limit where $N$ tends to infinity (compare for example \eqref{eq:tr2} and \eqref{eq:tr2N}). We also observed some instances of a property of factorisation, see for example \eqref{eq:tr22fac}. The factorisation is due to a phenomenon of concentration, with the effect that, as $N$ tends to infinity, and provided one scales the scalar product on $\u(N)$ correctly (which we did), the Wilson loop functionals, that is, the normalised traces of the random holonomies, become deterministic. The limit is thus a number depending on a loop, and this function is relatively simple, at least when one is working on the plane, because it satisfies, and is essentially determined, by the Makeenko--Migdal equations.

The main theorem of convergence is the following.

\begin{theorem}[Master field] \label{thm:CVMF} Let $M$ be either the plane $\R^{2}$ or the sphere $S^{2}$. For each $N\geq 1$, let $(H_{N,\ell})_{\ell\in \Loop(M)}$ be the Yang--Mills holonomy process on $M$ with structure group $G=\U(N)$, and with scalar product $\langle X,Y\rangle=N\Tr(X^{*}Y)$ on $\u(N)$. Then for every loop $\ell\in \Loop(M)$, the convergence of complex-valued random variables
\begin{equation}\label{eq:CVMF}
\tr(H_{N,\ell}) \build{\longrightarrow}_{N\to \infty}^{P} \Phi(\ell)
\end{equation}
holds in probability, towards a deterministic real limit.
\end{theorem}

This theorem was proved in \cite{LevyMF} in the case of the plane, and in \cite{DahlqvistNorris} in the case of the sphere, see also \cite{Hall}. In the case of the plane, which is simpler, it is also known that the convergence occurs quickly, in the sense that the series $\sum_{N\geq 1} {\rm Var}(\tr(H_{N,\ell}))$ converges. Thus, the convergence \eqref{eq:CVMF} holds almost surely. The conclusion is also known to be true if one replaces the unitary group by the special unitary group, the special orthogonal group, or the symplectic group. 

It is expected that Theorem \ref{thm:CVMF} is true on any compact surface, but a proof of this fact still has to be given.

In any case, when this theorem holds, the aforementioned asymptotic factorisation takes place, in the sense that for all loops $\ell_{1},\ldots,\ell_{n}$,
\[\lim_{N\to \infty}\E[\tr(H_{\ell_{1}})\ldots \tr(H_{\ell_{n}})]=\lim_{N\to \infty}\E[\tr(H_{\ell_{1}})]\ldots \lim_{N\to \infty}\E[\tr(H_{\ell_{n}})]=\Phi(\ell_{1})\ldots \Phi(\ell_{n})\]

The function $\Phi:\Loop(M)\to \R$ which appears in \eqref{eq:CVMF} is called the {\em master field}. This is a continuous function with respect to the convergence of loops with fixed endpoints (see the beginning of Section \ref{sec:contlim}) and it satisfies, crucially, the Makeenko--Migdal equation \eqref{eq:MMinf}, which is all that there is left of the full set of equations stated in Theorem \ref{thm:MME} as $N$ tends to infinity.

\begin{theorem} Assume that $M$ is either the plane $\R^{2}$ or the sphere $S^{2}$. The function $\Phi:\Loop(M)\to \R$ is the unique function that is continuous, invariant under area-preserving diffeomorphisms, satisfying the Makeenko--Migdal equation \eqref{eq:MMinf} and such that for every simple loop $\ell$ enclosing a domain of area $t$, one has, depending on whether $M$ is the plane or a sphere of total area $T$,
\begin{equation}\tag{$M=\R^{2}$}
\Phi(\ell)=e^{-\frac{t}{2}}
\end{equation}
or
\begin{equation}\tag{$M=S^{2}$}
\Phi(\ell)=\frac{1}{\pi} \int_{\R} \cosh \Big(\frac{x}{2}(T-2t)\Big) \sin (\pi \rho_{T}(x))\d x
\end{equation}
\end{theorem}

\subsection{A value of the master field on the plane}

As a conclusion to these notes, we give an example of computation of a value of the master field $\Phi$ on the plane, and choose an example that is not listed at the end of \cite{LevyMF}. We choose the loop $\ell$ represented on the left half of Fig. \ref{fig:exMF} below.

\begin{figure}[h!]
\begin{center}
\includegraphics{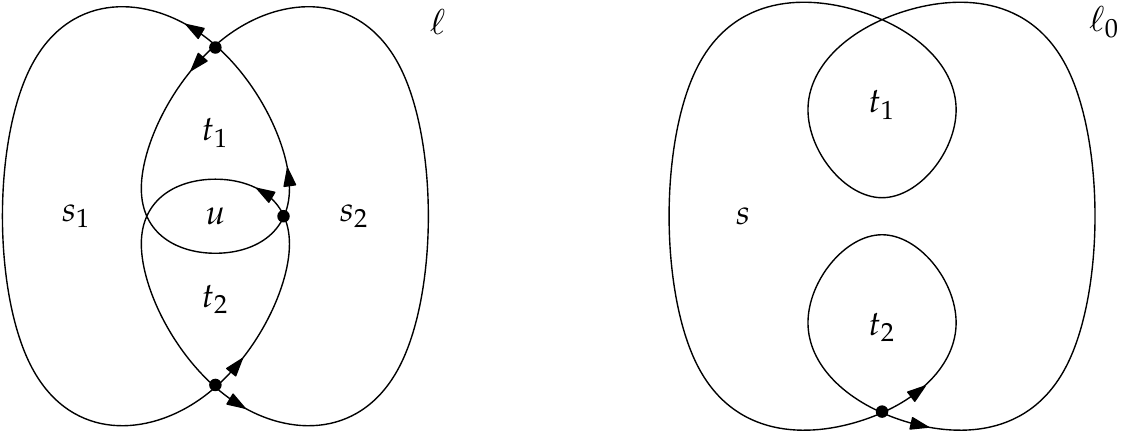}
\caption{\label{fig:exMF} \small We are interested in computing $\Phi(\ell)$. The strategy is to use the Makeenko--Migdal equations to compute $\partial_{u}\Phi(\ell)$. As $u=0$, the two inner windings of $\ell$ disentangle, and $\ell$ becomes identical to $\ell_{0}$. This loop $\ell_{0}$ is similar to the loop that we studied in Section \ref{sec:Yin}, and becomes exactly this loop when $t_{2}=0$. Our first task is thus to compute $\partial_{t_{2}} \Phi(\ell_{0})$.}
\end{center}
\end{figure}

Although we did not include this in our description of the function $\Phi$ on the plane $\R^{2}$, it is not difficult to check that the derivative of $\Phi$ of any loop with respect to the area of a face adjacent to the unbounded face is equal to $-\frac{1}{2}$ times the value of $\Phi$ on this loop. This factor $-\frac{1}{2}$ comes of course from the stochastic differential equation \eqref{eq:EDS} satisfied by the Brownian motion on $\U(N)$.

Given the value of $\Phi$ on simple loops and \eqref{eq:trheartMF}, the Makeenko--Migdal equation applied to the vertex of $\ell_{0}$ that is marked in Fig.~\ref{fig:exMF} yields
\[(2\partial_{s}-\partial_{t_{2}})\Phi(\ell_{0})=(-1-\partial_{t_{2}})\Phi(\ell_{0})=e^{-\frac{s}{2}-t_{1}-t_{2}}(1-t_{1})\]
which is solved in
\[\Phi(\ell_{0})=e^{-\frac{s}{2}-t_{1}-t_{2}}(1-t_{1})(1-t_{2})\]

If we can determine $\partial_{u}\Phi(\ell)$ explicitly, we are done, since $\Phi(\ell_{0})$ is exactly the value of $\Phi(\ell)$ at $u=0$. Applying the Makeenko--Migdal equations at the three marked vertices in Fig. \ref{fig:exMF} yields the derivatives $(\partial_{s_{1}}+\partial_{s_{2}}-\partial_{t_{2}})\Phi(\ell)$, $(\partial_{s_{1}}+\partial_{s_{2}}-\partial_{t_{1}})\Phi(\ell)$, and $(\partial_{t_{1}}+\partial_{t_{2}}-\partial_{s_{2}}-\partial_{u})\Phi(\ell)$.
Adding the three expressions and using the fact that $\partial_{s_{1}}\Phi(\ell)=\partial_{s_{2}}\Phi(\ell)=-\frac{1}{2}\Phi(\ell)$, we find
\[\Big(-\frac{3}{2}-\partial_{u}\Big)\Phi(\ell)=e^{-\frac{s_{1}+s_{2}}{2}-t_{1}-t_{2}-\frac{3u}{2}}(3-t_{1}-t_{2}-u)\]
and finally
\begin{empheq}[box=\mybox]{equation} \label{eq:exMF}
\Phi(\ell)=e^{-\frac{s_{1}+s_{2}}{2}-(t_{1}+t_{2})-\frac{3u}{2}}\Big(\frac{u^{2}}{2}+(t_{1}+t_{2}-3)u+(1-t_{1})(1-t_{2})\Big)
\end{empheq}
Evaluating this expression with $s_{1}=s_{2}=t_{1}=t_{2}=0$ yields the large $N$ limit of the third moment of the unitary Brownian motion at time $u$, as expressed by \eqref{eq:limmomMB} with $n=3$. This is consistent with the fact that shrinking all faces but the face of area $u$ reduces $\ell$ to a loop winding three times around a simple domain of area $u$.

\bibliographystyle{./Martin}
\bibliography{./bibSZ}

\end{document}